\documentclass[notitlepage,leqno]{article}

\usepackage[latin1]{inputenc}
\usepackage{latexsym}
\usepackage{amssymb}
\usepackage[T1]{fontenc}	
\usepackage[english]{babel}
\usepackage{amsmath,amssymb}
\usepackage{amsfonts}
\usepackage{amscd}
\usepackage{empheq}
\usepackage{fancyvrb}
\usepackage{marvosym}
\usepackage{mathrsfs}
\usepackage{amsthm}
\usepackage{hyperref}
\usepackage{verbatim}
\usepackage{tikz}
\usetikzlibrary{chains}
\usepackage{color}
\usepackage{chngpage}
\usetikzlibrary{matrix,through,arrows}

\renewcommand{\a }{\alpha }
\renewcommand{\b }{\beta }
\renewcommand{\d}{\delta }
\newcommand{\D }{\Delta }

\newcommand{\e }{\varepsilon }
\newcommand{\g }{\gamma}

\renewcommand{\l }{\lambda }

\newcommand{\s }{\sigma }
\newcommand{\Si }{\Sigma}

\newcommand{\be}{\begin{equation}}
\newcommand{\ee}{\end{equation}}
\newenvironment{pf}{\noindent{\sc Proof}.\enspace}{\rule{2mm}{2mm}\medskip}
\newenvironment{pfn}{\noindent{\sc Proof}}{\rule{2mm}{2mm}\medskip}
\newtheorem{thm}{Theorem}[section]

\newtheorem{lem}[thm]{Lemma}
\newtheorem{rem}[thm]{Remark}
\newtheorem{cor}[thm]{Corollary}

\newcommand{\R}{\mathbb{R}}
\newcommand{\Ss}{\mathbb{S}}

\renewcommand{\epsilon}{\varepsilon}

\newcommand{\nota}{\color{red}}
\newcommand{\fn}{\normalcolor}

\title{EXISTENCE OF  GENERALIZED TOTALLY UMBILIC 2-SPHERES IN PERTURBED 3-SPHERES}
\date{}

\author{Alessandro Carlotto\
   \thanks{\textsf{Stanford University, Department of Mathematics - Sloan Hall, 94305 Stanford, CA
}}
  \and
 Andrea Mondino\
 \thanks{\textsf{ETH Zentrum, Department of Mathematics,  CH-8093 Z\"urich, Switzerland}}
  }
\begin{document}

\maketitle

\

\begin{description}
\item[Abstract.]
It was recently shown by R. Souam and E. Toubiana \cite{ST} that the (non constantly curved) Berger spheres do not contain totally umbilic surfaces. Nevertheless in this article we show, by perturbative arguments, that all analytic metrics sufficiently close to the round metric $g_{0}$ on $\mathbb{S}^{3}$ possess \textsl{generalized} totally umbilic 2-spheres, namely critical points of the conformal Willmore functional $\int_{\Si}|A^\circ|^{2}\,d\mu_{\gamma}$. The same is true in the smooth setting provided a suitable non-degeneracy condition on the traceless Ricci tensor holds. The proof involves a gluing process of two different finite-dimensional reduction schemes, a sharp asymptotic analysis of the functional on perturbed umbilic spheres of small radius and a quantitative Schur-type Lemma in order to treat the cases when the traceless Ricci tensor of the perturbation is degenerate but not identically zero. For left-invariant metrics on $SU(2)\cong \mathbb{S}^{3}$ our result implies the existence of uncountably many distinct Willmore spheres. 
\end{description}

\

\section{Introduction}
One of the most general problems in \textsl{extrinsic} Riemannian geometry is to find the \textsl{best immersion} $\Phi$ of a given smooth manifold $M$ in a higher dimensional ambient space $N$ endowed with a Riemannian metric $g$, by which we mean that some special curvature condition is required. For instance, one might prescribe the vanishing of the second fundamental form $A$ of $M$ in $N$ (in which case $\Phi$ is called \textsl{totally geodesic}), or the vanishing of its trace $H$ (in this case $\Phi$ is \textsl{minimal}) or, instead, the vanishing of its traceless component $A^\circ$, which corresponds to $\Phi$ being \textsl{totally umbilic}. Over the last century, a number of different obstructions to the existence of such \textsl{optimal} immersions have been found: for example, it is well-known that if $M$ is compact and has no boundary then it cannot be minimally immersed in the Euclidean space $\R^{n}$, as a basic consequence of the classical monotonicity formula for minimal submanifolds.
If $\textrm{dim}_{\R}\left(M\right)=2$ and $\textrm{dim}_{\R}\left(N\right)=3,$ a recent  obstruction was found by R. Souam and E. Toubiana \cite{ST}, who proved that if $(N,g)$ belongs to the class of the \textsl{Berger spheres} then there exist no totally umbilic immersions of $M$ in $N$, unless $(N,g)$ is in fact a space form. Given this fact, it is then natural to weaken our requirement and to ask whether there exist immersions having the property to minimize (or, more generally, to be critical points for) some integral functional of $|A^\circ|$. This is the object of study of the present article.
\newline

Before proceeding further, let us introduce some notation. We will always deal with a compact, isometrically immersed, surface $(\Si,\g)$ in $(N,g)$: the corresponding principal curvatures will be denoted by $\l_{1}, \l_{2}$ and the mean curvature $H$ of $\Sigma$ will be their \textsl{sum}, namely $H=\l_{1}+\l_{2}$. Moreover, according to our sign convention the round unit sphere in $\R^{3}$ has mean curvature equal to $2$.  \newline
 
In this work, we consider the \textsl{conformal Willmore functional} given by
\be\label{def:I}
I\left(\Si,\g\right)=\int_{\Si}\left(\frac{H^{2}}{4}-D\right)\,d\mu_{\g}=\frac{1}{2}\int_{\Si}|A^\circ|^{2}\,d\mu_{\g}
\ee
where $d\mu_{\g}$ is the Riemannian volume form associated to the metric $\g$, $A^\circ:=A-\frac{1}{2} H \g$ is the traceless part of the second fundamental form $A$ of $\Sigma$ in $M$ and $D=\l_{1}\l_{2}$. 
In the special case when the sphere $\mathbb{S}^{3}$ is endowed with its standard round metric $g_{0}$ (as embedded unit sphere in $\R^{4}$) then the previous two functionals coincide, modulo a null Lagrangian, with the functional
\be
W\left(\Sigma,\g\right)=\int_{\Si}\left(\frac{H^{2}}{4}+1\right)\,d\mu_{\g}
\ee
as an immediate consequence of the Gauss equations. However, from the point of view of conformal geometry, for immersions in a general Riemannian manifold, the functional $I$ is  more natural than $W$, indeed the former is \textsl{conformally invariant} (see for instance \cite{Wei}) while the latter may not be.
\\ Since we deal with isometric immersions, $\g$ is nothing but the pullback of the metric $g$,  hence we will omit the $\g$-dependence of $I$ and $W$ (therefore we will only write $I(\Sigma)$ and $W(\Sigma)$). 

\

In the last five decades, the study of the existence of critical points for these functionals and their geometric characterization has been the object of a number of works also due to the connections with other  key questions in Geometric Analysis, like the classification of positive genus minimal surfaces in the round 3-sphere $\mathbb{S}^{3}$ (this link, inspired by the works of A. Ros \cite{Ros} and F. Urbano \cite{Urb}, was smartly exploited by F. Marques and A. Neves in their proof of the Willmore conjecture  \cite{MN}), the  regularity of complete properly embedded minimal surfaces in the Hyperbolic $3$-space (see the recent paper of S. Alexakis and R. Mazzeo \cite{AlMa}),  the study of sharp eigenvalue estimates in relation to the conformal volume theory by P. Li and S. T. Yau (\cite{LY}), etc.

Starting with the nowadays classical paper of  L. Simon \cite{SiL}, the variational study of Willmore-type functionals for immersions \emph{in the flat  Euclidean  space} $\R^n$ has been extensively carried through, both in terms of existence and of regularity results  (we recall the remarkable results of E. Kuwert and R. Sch\"atzle, for instance \cite{KS}). In parallel with Simon's ambient approach (involving geometric measure theory), recently T. Rivi\`ere  developed a parametric approach to successfully attach the existence and regularity issues regarding the Willmore functional (see \cite{Riv} and \cite{Riv2}).

We stress that all the aforementioned existence results concern  immersions \emph{in the flat euclidean space $\R^n$}  (or equivalently, thanks to the stereographic projections and the conformal invariance of the Willmore functional, for immersions in the round sphere $\Ss^n$). Explicit examples of Willmore surfaces, or explicit bounds on the energy of some special submanifolds,  in very symmetric ambient manifolds have been constructed by several authors (see for instance  \cite{BaFe}, \cite{HS}, \cite{LU},\cite{Wang}, etc.). The existence of Willmore surfaces in non constantly curved ambient manifolds is a very recent topic started by the second author in a perturbative setting  in \cite{Mon1} and \cite{Mon2} (for existence of Willmore surfaces under area constraint, still in a perturbative setting, see the papers of T. Lamm, J. Metzger and F. Schulze \cite{LMS}, \cite{LM} and \cite{LM2}).  The minimization (i.e. the existence of a minimizer and the corresponding regularity theory), among smooth immersions of $2$-spheres in a $3$-manifold, of quadratic curvature functionals of the type $\int |A|^2$ and $\int |H|^2+1$ has been achieved by the second author in collaboration with E. Kuwert and J. Schygulla  in \cite{KMS} (see also \cite{MonSchy} for the non compact case).  Finally, in collaboration with T. Rivi\`ere (see \cite{MoRi1} and \cite{MoRi2}), the second author developed a parametric approach for studying the regularity of possibly branched immersions which are critical points (possibly with constraints) of such curvature functionals in Riemannian manifolds (also of higher codimension), and applied this theory to  the minimization of these functionals among possibly branched immersions of $2$-spheres in homotopy groups (see also the paper of J. Chen and Y. Li \cite{ChenLi} for related results). 

\

Let us remark that all these results in general Riemannian manifolds follow a \emph{minimization} scheme; on the other hand such a method cannot be applied for finding interesting spherical type critical points of the conformal Willmore functional $I$ defined in \eqref{def:I}. Indeed,  given a point $p$ in the $3$-manifold $M$ and denoted with $S_{p,\rho}$ the geodesic sphere of center $p$ (i.e. the sphere in geodesic coordinates centered at $p$)  and radius $\rho>0$, it is easy to see that $I(S_{p,\rho})\to 0$ as $\rho \to 0$. Therefore the infimum, among smooth immersions of $\Ss^2$ into $M$, of the functional $I$ is zero and every minimizing sequence either collapses to a point or converges, in a suitable sense, to a totally umbilic surface. But both the situations are not interesting for our purposes: the former is a degeneration which does not give a geometric object in the limit, while the latter simply may not happen (and indeed it \textsl{does not} happen in our case, by the result of Souam and Toubiana).
 
Therefore, in order to study the critical points of the functional $I$,  one could either  perform a $\min-\max$ scheme or  use a perturbative method. The present paper is related to the second technique,  the first one will be studied in a forthcoming work.  

Now we can state the main theorem of this paper which answers, in a perturbative setting, to the question of the existence of an \textsl{umbilically best immersion} in relation to the aforementioned obstruction given in \cite{ST}.

\begin{thm}
\label{main}
Let $g_\e=g_{0}+\e h$ be a Riemannian metric on $\mathbb{S}^{3}$ for some analytic, symmetric $(0,2)$-tensor $h$. There exists $\overline{\e}\in\mathbb{R}_{>0}$  such that if $\e\in\left(-\overline{\e},\overline{\e}\right)$ then there exist embedded critical points for the conformal Willmore functional in metric $g_\e$
\be\nonumber
I_{\e}\left(\Si\right)=\frac{1}{2}\int_{\Si}|A^\circ|^{2}\,d\mu_{\g_{\e}}.
\ee
More precisely, every Willmore surface we construct is a normal graph over a totally umbilic sphere of a smooth function $w_\e$ converging to $0$ in ${C^{4,\alpha}}$ norm as $\e \to 0$. 
\end{thm}

\begin{rem}
It follows from our construction (as explained in the sequel of this Introduction) that the critical points for which we show existence are in fact \textsl{saddle points} for $I_{\e}$. Moreover, a standard \textsl{bumpy-metric argument} shows that (in case $(\mathbb{S}^{3}, g_{\epsilon})$ \emph{does not}  have constant sectional  curvature) these are \textsl{generically non-degenerate} of index exactly 4. To our knowledge, this is the first existence result for embedded Willmore surfaces of saddle type in a compact (non constantly curved) ambient manifold.
\end{rem}

Let us now briefly describe the logical scheme of the proof and, correspondingly, the structure of this article. The first basic idea is, loosely speaking, to compare the functional $I_{\e}$ to its unperturbed counterpart $I_{0}$: indeed, if $\e=0$ then we endow $\mathbb{S}^{3}$ with the round metric $g_{0}$ and clearly the conformal Willmore functional admits a four-dimensional manifold of minimum points (where $I_{0}$ is identically null) which is made of all totally umbilical spheres $(S_{p,\rho})$ for $p\in \mathbb{S}^{3}$ and $\rho\in\left(0,\pi\right)$. In fact, thanks to the quantization results by R. Bryant \cite{bry} (see also the recent paper  \cite{LaNg} of T. Lamm and H. T. Nguyen for the branched case) it is well-known that we can \textsl{separate such manifold}: the full moduli space of smooth immersions of $\Ss^2$ into $\Ss^3$ which are critical points  for the functional $W_{0}$  (or, equivalently, to $I_0$) consists of a countable number of connected components, the one of minimal $W_{0}-$energy corresponding to the \textsl{totally umbilic} 2-spheres in $\mathbb{S}^{3}$ (for which $W_{0}=4\pi$), followed by a second component with a gap of exactly $12\pi$. As a result, for our purposes we can neglect all higher energy components and so, from now onwards, let us denote by $Z'$ the (closure of the) critical sub-manifold of $I_{0}$ where such functional vanishes identically. The perturbation scheme we need to apply is based on the explicit knowledge of the global topology of the critical manifold $Z'$:
unfortunately, if $\Sigma$ is \textsl{not oriented}, then we \textsl{cannot} simply identify $Z'$ with $\mathbb{S}^{3}\times \left[0,\pi\right]$ since clearly the couples $(p,r)$ and $(\hat{p},\pi-r)$ correspond to \textsl{the same sphere} (where we have denoted by $\hat{p}$ the antipodal point of $p$ in $\mathbb{S}^{3}$).
Instead, the set $Z'$ is diffeomorphic to $\textrm{Bl}_{O}D^{4}$, namely the \textsl{real blow-up} at the origin of the unit disk $D^{4}\subseteq\mathbb{R}^{4}$. We recall that
\be\nonumber
\textrm{Bl}_{O}D^{4}=\left\{\left(x,l\right)\in D^{4}\times\mathbb{R}\mathbb{P}^{3} | \ x\in l\right\}
\ee
and that $\textrm{Bl}_{O}D^{4}\cong D^{4}\# \mathbb{R}\mathbb{P}^{4}$  so that, as a result, $Z'$ is a closed smooth \emph{non-orientable} manifold with boundary.
\newline 
Given this fact and in order to avoid unnecessary complications in our proof, we will therefore consider the same Willmore functional $I_{\e}$ defined on \textsl{oriented} isometrically embedded spheres, so that the critical manifold at minimal energy for $W_{0}$ (which we denote by $Z$) is diffeomorphic to $\mathbb{S}^{3}\times \left[0,\pi\right]$. We remark that this choice is not at all necessary for our arguments to work, yet in the former setting the logical structure of our proof would considerably lose in terms of conceptual clarity and effectiveness.
For the sake of brevity, we will sometimes refer at $Z$ as a \textsl{cylinder} and its subsets $\mathbb{S}^{3}\times\left\{0\right\}$ and $\mathbb{S}^{3}\times \left\{\pi\right\}$ will be called \textsl{bases} of such cylinder.

\

Given such manifold $Z$, our strategy is based on applying a finite-dimensional reduction of our problem in the spirit of Ambrosetti-Badiale \cite{ab1}-\cite{ab2} (which in turn is based on the classical Lyapunov-Schmidt reduction; for a sketch of the abstract method see Subsection \ref{LS}, for a more extensive discussion including applications see \cite{am2}): namely we construct a small perturbation $Z_{\e}$ of $Z$ which plays the role, for the functional $I_{\e}$, of a \textsl{natural variational constraint} in the sense that (interior) critical points of $I_{\e}$ on $Z_{\e}$ are in fact critical points for the unconstrained functional $I_{\e}$. The manifold $Z_{\e}$ corresponds, in our setting, to small \textsl{graphical} perturbations of the totally umbilic spheres $S_{p,\rho}$: we construct a map $w_{\e}(p,\rho)$ that associates, to each point $p \in\mathbb{S}^{3}$ and radius $\rho$, a (smooth) function defined over the 2-sphere $\mathbb{S}^{2}$, identified with $S_{p,\rho}$. In fact, in performing this reduction we face a number of technical obstacles. First of all, the general method by Ambrosetti-Badiale refers to a Hilbert space setting, while in our case it is convenient to work with $w\in C^{4,\a}\left(\mathbb{S}^{2};\R\right)$ (the equation we need to solve in order to find critical points of $I_{\e}$ has order four) and therefore we need to adapt the construction to our specific setting. The second issue (which is related to the first) is that the construction of the map $w_{\e}(\cdot,\cdot)$ (see Section \ref{fdred}) requires the second derivative operator of $I_{\e}$ to be \textsl{uniformly elliptic}, which is true (by equation \eqref{secvar}) only on relatively compact subsets of $\mathring{Z}\simeq \mathbb{S}^{3}\times\left(0,\pi\right)$. As a result, in order to solve our problem by means of a finite dimensional reduction we first need to show that \textsl{we can get rid of suitably small neighborhoods of the bases of $Z$}. More precisely, we will show that the functional $I_{\e}$ is strictly increasing for $\rho\simeq 0$ (hence, by symmetry, for $\rho\simeq\pi$) so that we can apply the reduction scheme to a suitable closed subcylinder of $Z$ of the form $\mathbb{S}^{3}\times \left[\d, \pi-\d\right]$ for some suitably small $\d>0$. A delicate aspect, in doing this, is that we need to find $\d$ \textsl{not depending on $\e$}. To this aim we need to construct (in Section \ref{srasy}) a sort of \textsl{second} finite-dimensional reduction map (which we will still denote by $w_{\e}(\cdot,\cdot)$) in order to study the asymptotics for $\rho\to 0$ of $I_{\e}$ on graphs over the spheres $S_{p,\rho}$ for very small values of their radii. In other terms, we show that if $I_{\e}$ has a critical point $\tilde{w}$ which is a graphical perturbation of a totally umbilic sphere $S_{p,\rho}$ of small radius then in fact $\tilde{w}=w_{\e}\left(p,\rho\right)$ and at that point we study the behavior, both in $\rho$ and in $\e$, of the corresponding reduced functional $I_{\e}\left(S_{p,\rho}\left(w_{\e}\left(p,\rho\right)\right)\right)$. This is based on the work performed in Section 3 of \cite{Mon2} (where it is proved that $I_{\e}\left(S_{p,\rho}\left(w_{\e}\left(p,\rho\right)\right)\right)\simeq \frac{\pi}{5}\rho^{4}|\mathring{\textrm{Ric}}|^{2}+O_{\e}(\rho^{5})$); yet in our setting we do not fix a given Riemannian metric on $\mathbb{S}^{3}$ and it is crucial for us to obtain estimates that are uniform in $\e$, at least for small values of this parameter.

\

A distinctive feature of our Theorem \ref{main}, compared to other perturbative results (like, for instance, \cite{Mon1} and \cite{Mon2}), is that the perturbation $h$ is \textsl{completely arbitrary} in the sense that we do not add any sort of technical non-degeneracy condition.
This is possible thanks to the fact that if the traceless Ricci tensor of $g_{\e}$ vanishes identically on $\mathbb{S}^{3}$, then $g_{\e}$ is homothetic to the round metric;  if this is not the case, we can exploit the fact that all curvature tensors of $g_{\e}$ are in fact analytic in $\e$ and hence $\mathring{\textrm{Ric}}$ can be expanded in $\e$ (with analytic coefficients) and the problem can be suitably reduced to the fully degenerate case thanks to a quantitative Schur-type argument (we remark that an \emph{integral}-quantitative Schur Lemma was proven by C. De Lellis and P. Topping \cite{DLT} in case of positive Ricci curvature, but our arguments are independent from theirs: indeed on one hand we work in a perturbative regime, on the other hand we get \emph{pointwise} estimates). We believe that the method used in this part of the argument is rather new and interesting in itself. 

\

When the analytic assumption on $g_{\e}$ (or, equivalently, on $h$) is removed, our argument still works provided we require that the expansion (in $\e$) of the traceless Ricci tensor in \textsl{non degenerate} at least at some point. This amounts to requiring that the variation $h$ is not in the kernel of the linearization at the metric $g_{0}$ of the traceless Ricci operator. Therefore, we can state the following smooth counterpart of our main result.

\begin{thm}
\label{smooth}
Let $g_\e=g_{0}+\e h$ be a Riemannian metric on $\mathbb{S}^{3}$ for some $C^{\infty}$ symmetric $(0,2)$-tensor $h$. There exists $\overline{\e}\in\mathbb{R}_{>0}$ such that if $\e\in\left(-\overline{\e},\overline{\e}\right)$ then there exist embedded critical points for the conformal Willmore functional in metric $g_\e$
provided $h$ does \textsl{not} satisfy the following equation
\be\nonumber
\frac{1}{2}\D_{L}h+\frac{1}{2}\mathcal{L}_{{\d G\left(h\right)}^{\sharp}}g_{0}-\frac{1}{3}g_{0}\left(Ric_{g_{0}},h\right)+\frac{1}{3}\left(\d^{2}h\right)g_{0}-\frac{1}{3}\D(\textrm{tr}_{g_{0}}h)g_{0}=0.
\ee
Here $\D_{L}$ denotes the Lichnerowicz Laplace operator, $\d$ is the divergence (with respect to the metric $g_{0}$), $^{\sharp}$ is the standard musical isomorphism $\Gamma\left(T^{\ast}M\right)\to\Gamma\left(TM\right)$ determined by the metric $g_{0}$ and we have set $G(h)=h-\frac{1}{2}\left(\textrm{tr}_{g_{0}}h\right)g_{0}$.
\\
Moreover, every Willmore surface we construct is a normal graph over a totally umbilic sphere of a smooth function $w_\e$  converging to $0$ in ${C^{4,\alpha}}$ norm as  $\e \to 0$. 
\end{thm}
\

Our argument also implies multiplicity results whenever $\mathbb{S}^{3}$ is endowed with left-invariant metrics with respect to its Lie group structure. Indeed, the manifold $\mathbb{S}^{3}$ can be identified with the algebraic group 
\be
\nonumber
SU(2)=\left\{A\in M_{2\times 2}\left(\mathbb{C}\right):\ \det(A)=1, \ A^{\ast}=A^{-1}\right\}
\ee
\be
\nonumber
=\left\{\begin{bmatrix}
 z & w \\
 -\overline{w} & \overline{z}
\end{bmatrix}
: \ \left|z\right|^{2}+\left|w\right|^{2}=1\right\}.
\ee
Now, if $G$ is an analytic Lie group (in fact we know that any $C^{0}$ Lie group always admits a unique analytic structure) then every left-invariant Riemannian metric is itself analytic and therefore we are always in position to apply our Theorem \ref{main} to left-invariant metrics on $SU(2)$ without any non-degeneracy constraint. If we combine this fact with the trivial remark that the action of $G=SU(2)$ on itself is \textsl{transitive} we obtain the following remarkable consequence.

\begin{cor}
\label{group}
Let $g_{\e}=g_{0}+\e h$ be a left-invariant metric on $SU(2)\cong S^{3}$. There exists $\overline{\e}\in\mathbb{R}_{>0}$ such that if $\e\in\left(-\overline{\e},\overline{\e}\right)$ then for every $p\in\mathbb{S}^{3}$ there exists an embedded critical 2-sphere for the conformal Willmore functional (in metric $g_{\e}$) passing through $p$. As a result, under these assumptions the functional $I_{\e}$ has uncountably many distinct critical points.
\end{cor}

\begin{rem}
\label{expas}
We would like to stress that both Theorem \ref{main} and Corollary \ref{group} can easily be extended, with almost no change in the proofs, to the case when $g_{\e}=g_{0}+h_{\e}$ with $h$ a perturbation which is analytic in all of its variables and such that $h_{\e=0}=0$. Therefore, it is clear that Corollary \ref{group} does in fact apply to \textsl{any} left-invariant metric $g$ on $SU(2)$ which is sufficiently close to $g_{0}$. 
\end{rem}

\begin{rem}
It is also appropriate to remark that the Berger spheres are indeed a one parameter family of left-invariant metrics on $\mathbb{S}^{3}\cong SU(2)$, so that Corollary \ref{group} implies the existence for them of uncountably many \textsl{generalized} totally umbilic 2-spheres, in sharp contrast with the negative result of Souam and Toubiana asserting that there are no totally umbilic 2-spheres at all.
More generally, our multiplicity result apply to the subclass of left-invariant metrics $\mathcal{G}$ defined by requiring that
\be\nonumber
g\in \mathcal{G} \ \Rightarrow \ g(X_{i}, X_{j})=\d^{i}_{j}\l_{i} \ \textrm{for some} \ \l_{1},\l_{2},\l_{3}\in\mathbb{R}_{>0}
\ee
where 
\be\nonumber
X_{1}=\begin{bmatrix}
 i & 0 \\
 0 & -i
\end{bmatrix}, 
X_{2}=\begin{bmatrix}
 0 & 1 \\
 -1 & 0
\end{bmatrix},  
X_{3}=\begin{bmatrix}
 0 & i \\
 i & 0
\end{bmatrix}.
\ee
These clearly form a basis of the Lie algebra $\mathfrak{su}\left(2\right)$ of $SU(2)$. Notice that the Berger spheres correspond to the 1-parameter family in $\mathcal{G}$ given by choosing $\l_{1}=\l$, $\l_{2}=\l_{3}=1$ and in that case $X_{1}$ is tangent to the orbits of the Hopf circle action.
\end{rem}

\
The paper is structured as follows: in Section 2 we collect some preliminary results (both concerning perturbation schemes and expansions of curvature tensors), in Section 3 we construct the manifold $Z_{\e}$, namely the finite-dimensional reduction map $w$, in Section 4 we study the sharp asymptotics of the functional $I_{\e}$ for small radii and finally we give in Section 5 a detailed proof of Theorem \ref{main}. The first and second variation formulas for the conformal Willmore functional, which are recalled in Section 2, are proved in the Appendix at the end of this work.

\

\textsl{Acknowledgments}.  During the preparation of this work, A. C. was supported by  NSF grant DMS/0604960 and A. M. was  supported by the ETH-Fellowship. The authors would like to thank Andrea Malchiodi for introducing them to variational and perturbative methods in Nonlinear Analysis. Gratitude is also expressed to Otis Chodosh for his careful reading of the preliminary manuscript.

\section{Notation and preliminary results}

In order to make more concise and readable the key arguments in the  proof of Theorem \ref{main}, we collect in this section a number of useful results.  For each of them, we will either provide a proof or give the reader an appropriate reference.

\subsection{Notations}
\label{nota}
As anticipated in the Introduction, it is convenient for us to consider small perturbations of the totally umbilic spheres $S_{p,\rho}$ in $\mathbb{S}^{3}$ of center $p\in \Ss^3$ and radius $\rho\in[0,\pi]$; such perturbations are of the form of normal graphs defined over the unit sphere $\mathbb{S}^{2}\hookrightarrow \R^{3}$ (the identification using the exponential map of the appropriate metric $g_{\e}$). Coherently with \cite{Mon1} and \cite{Mon2}, we will denote by $\Theta_{1}$ and $\Theta_{2}$ the corresponding coordinate vector-fields on $\mathbb{S}^{2}$  (induced by the standard polar coordinates on the unit sphere of the Euclidean space $\R^{3}$).
Due to technical reasons (specifically: the need to apply suitable Schauder estimates), we will take $w\in C^{4,\a}\left(\mathbb{S}^{2};\R\right)$, which is the Banach space of functions whose $4^{th}$ order derivatives with respect to $\Theta_{i}, \ i=1,2$, are $\a-$H\"older, for some specific $\a\in \left(0,1\right)$. Denoted by $\D_{\mathbb{S}^{2}}$ the Laplace-Beltrami operator on $\mathbb{S}^{2}$, we will often work with the fourth order operator $\D_{\mathbb{S}^{2}}\left(\D_{\mathbb{S}^{2}}+2\right)$, which induces a splitting of the Hilbert space $L^{2}\left(\mathbb{S}^{2}\right)$ as follows:
\be\nonumber
L^{2}\left(\mathbb{S}^{2}\right)=\textrm{Ker}\left[\D_{\mathbb{S}^{2}}\left(\D_{\mathbb{S}^{2}}+2\right)\right]\oplus\textrm{Ker}\left[\D_{\mathbb{S}^{2}}\left(\D_{\mathbb{S}^{2}}+2\right)\right]^{\perp}.
\ee 
We shall then consider $C^{4,\a}\left(\mathbb{S}^{2}\right)$ as a subspace of $L^{2}\left(\mathbb{S}^{2}\right)$, hence there is an induced splitting as above and we can set
\be
\label{split}
C^{4,\a}\left(\mathbb{S}^{2}\right)^{\perp}=C^{4,\a}\left(\mathbb{S}^{2}\right)\cap \textrm{Ker} \left[\D_{\mathbb{S}^{2}}\left(\D_{\mathbb{S}^{2}}+2\right)\right]^{\perp}.
\ee
We remark that $K=\textrm{Ker}\left[\D_{\mathbb{S}^{2}}\left(\D_{\mathbb{S}^{2}}+2\right)\right]$ is finite dimensional therefore closed, and $C^{4,\a}\left(\mathbb{S}^{2}\right)^{\perp}$ is itself a Banach space with (the restriction of) the $C^{4,\a}-$norm.
Finally, it is also convenient to name $P:L^{2}\left(\mathbb{S}^{2};\R\right)\to K^{\perp}$ the $L^{2}-$orthogonal projector to $K^{\perp}.$

\

Given a point $p\in \mathbb{S}^{3}$, $\rho\in\left(0,\pi\right)$ and a function $w\in C^{4,\a}\left(\mathbb{S}^{2};\R\right)$ (of suitably small norm), we are then in position to define a \textsl{perturbed geodesic sphere}, denoted by $S_{p,\rho}\left(w\right)$ as the image of the map $\Psi_{p,\rho,w,\e}:\mathbb{S}^{2}\to \mathbb{S}^{3}$
\be
\label{graph}
\Psi_{p,\rho,w,\e}(\Theta)=exp_{p}\left(\left(\rho+w\left(\Theta\right)\right)\Theta\right)
\ee
which is in fact a normal graph over $S_{p,\rho}\hookrightarrow \mathbb{S}^{3}$. We stress that here $exp_{p}$ denotes the exponential map defined on $T_{p}\mathbb{S}^{3}$ for a given metric $g_{\e}=g+\e h$ on $\mathbb{S}^{3}$, not necessarily the round one (in which case, we will add an explicit remark to our discussion). 

\

Given $a\in \mathbb{N}$, any expression of the form $L_{p}^{\left(a\right)}\left(w\right)$ denotes a \textsl{linear} combination of $w$ and its derivatives (with respect to $\Theta_{1}$ and $\Theta_{2}$) up to order $a$. We allow the coefficients of such combination to depend (smoothly) on $p, \rho$ and $\e$, but we require the existence of a constant $C$ (independent of these) so that 
\be\nonumber
\left\|L_{p}^{\left(a\right)}\left(w\right)\right\|_{C^{k,\a}\left(\mathbb{S}^{2}\right)}\leq C\left\|w\right\|_{C^{k+a,\a}\left(\mathbb{S}^{2}\right)}, \ k\in\mathbb{N}.
\ee
More generally, for $b\in\mathbb{N}$, any expression of the form $Q_{p}^{\left(b\right)\left(a\right)}\left(w\right)$ denotes a polynomial expression involving monomials of degree \textsl{at least} $b$, each of these involving $w$ and its derivatives up to order $a$. Again, we allow the coefficients to depend on our parameters, yet we require the existence of absolute constants giving bounds as above and also of the form
\be\nonumber
\left\|Q_{p}^{\left(b\right)\left(a\right)}\left(w_{2}\right)-Q_{p}^{\left(b\right)\left(a\right)}\left(w_{1}\right)\right\|_{C^{k,\a}\left(\mathbb{S}^{2}\right)}\leq C\left(\left\|w_{2}\right\|_{C^{k+a,\a}\left(\mathbb{S}^{2}\right)}+\left\|w_{1}\right\|_{C^{k+a,\a}\left(\mathbb{S}^{2}\right)}\right)^{b-1}\left\|w_{2}-w_{1}\right\|_{C^{k+a,\a}\left(\mathbb{S}^{2}\right)}
\ee   
provided $\left\|w_{l}\right\|_{C^{a}\left(\mathbb{S}^{2}\right)}\leq 1, \ l=1,2$. If the numbers $a, b$ are not specified we agree that  they equal 4 and 2 respectively.

\

If $x$ is a real variable and $f:I\to\R$ is a function of $x$ defined at least on some neighborhood of zero, we will write $f(x)=O(\left|x\right|^{\b})$ (for some $\b\in\R_{>0}$) in order to mean that 
$$\limsup_{x\to 0}\frac{\left|f(x)\right|}{|x|^\beta}< \infty.$$
When $f$ depends (smoothly enough) on some other variable, say $z$, we will use the notation $O_{z}\left(\left|x\right|^{\b}\right)$ in order to stress the dependence on $z$ (and, more specifically, to stress the fact that the remainder might not be uniform in $z$). In our problem, we need to consider functionals and functions depending on several parameters, typically $p\in \mathbb{S}^{3}$, $\rho\in\left[0,\pi\right]$ (or possibly in a smaller interval) and $\e\in\left(-\e^{*},\e^{*}\right)$ for some suitably small $\e^{*}$ and therefore we will often write $O_{p,\rho}\left(\e^{\b}\right)$ and $O_{p,\e}\left(\rho^{\b}\right)$ whenever an estimate is gotten by \textsl{freezing} some of the parameters (e.g. $p,\rho$ and $p,\e$ respectively) and considering the asymptotics with respect to the other ones.

\

\subsection{Riemannian geometry preliminaries}
\label{RGP}

Given a Riemannian metric $g$ on $\mathbb{S}^{3}$, we will only make use of the associated Levi-Civita connection $\nabla$ and concerning all the corresponding curvature tensors we will follow the conventions given, for instance, on the book by Petersen \cite{Petersen}. It is a trivial, yet crucial remark that if $g=g_{0}+\e h$ then all curvature tensors are \textsl{analytic} in $\e$: thus for $|\e|<\overline{\e}$ they can be expanded in power series of $\e$ with smooth coefficients if $h$ is, or more generally of class $C^{k}$ if $h$ is a tensor of class $C^{k+2}$. Both these statements follow at once from the local expression of the curvature tensors. 
In this work, we will mostly be interested in the Ricci curvature tensor $\textrm{Ric}_{g_\e}$ of $g_{\e}$ and in its trace-free part $\mathring{\textrm{Ric}}_{g_\e}:=\textrm{Ric}_{g_\e}-\frac{1}{3} R_{g_\e} g_\e$ where $R_{g_\e}$ is the scalar  curvature of the same metric. Concerning the perturbative expansion of $\mathring{\textrm{Ric}}_{g_\e}$, observe that

\be \nonumber
|\mathring{\textrm{Ric}}_{g_\e}|^{2}=\e^{2}T^{\left(2\right)}_{p}(h)+o(\e^{2})
\ee
where $T_{p}^{(2)}(h)$ denotes a non-negative quadratic expression in the second derivatives of $h$ and namely (see the statement of Theorem \ref{smooth})
\be \nonumber
T^{\left(2\right)}_{p}(h)=\left(\frac{1}{2}\D_{L}h+\frac{1}{2}\mathcal{L}_{{\d G\left(h\right)}^{\sharp}}g_{0}-\frac{1}{3}g_{0}\left(Ric_{g_{0}},h\right)+\frac{1}{3}\left(\d^{2}h\right)g_{0}-\frac{1}{3}\D(\textrm{tr}_{g_{0}}h)g_{0}\right)^{2}
\ee
Moreover if $T_{\cdot}^{(2)}(h)\equiv 0$ identically on $\mathbb{S}^{3}$ (which is a non-generic condition on the perturbation $h$) then locally (around any given point) $|\mathring{\textrm{Ric}}_{g_\e}|^{2}=\sum_{k\geq k_{0}}\e^{k}T^{\left(k\right)}_{p}(h)$ for some $k_{0}\geq 4$ and with suitably strong convergence in a (possibly smaller) neighborhood.

\subsection{First and second variation formulas}
\label{fs}

Given $p\in\mathbb{S}^{3}$ and $\rho\in(0,\pi)$ we state here the first and second variations of the functional $I_{0}$ on the totally umbilic spheres $S_{p,\rho}$ (with the pullback metric $\g_{0}$ given by the restriction of $g_{0}$) , the proof being postponed to Appendix A. \newline
\begin{lem}
\label{fslem}
Let us consider an isometrically immersed surface $(\Si,\g)$ and a deformation $F:\Si\times(-\s,\s)\to\mathbb{S}^{3}$ such that $F(\Si,0)=\Si$ and $\frac{\partial F}{\partial s}\left(\Sigma,0\right)=u\nu$ where $\nu$ is the (co-)normal vector field of $\Si$ (which is oriented, by assumption) in $\mathbb{S}^{3}$ and $u\in C^{4,\alpha}\left(\Sigma, \R\right)$. Then, if we set $L$ to be the Jacobi operator of $\Si$, namely
\be\nonumber
Lu=-\D_{\Si,\g}u-\left(Ric\left(\nu,\nu\right)+\left|A\right|^{2}\right)u \ ,
\ee we have that:
\begin{enumerate}
\item{the first variation formula for $I_{0}$ is given by
\be\label{first}
\d I_{0}\left(\Si\right)\left[u\right]=\int_{\Si}u\left[\frac{1}{2} LH+\left(\frac{H^{3}}{4}+H\right)\right]\,d\mu_{\gamma},
\ee
hence the first derivative operator is $I_{0}' \left(\Sigma\right)= \frac{1}{2} LH+\left(\frac{H^{3}}{4}+H\right)$, and $\left(\Sigma,\gamma\right)$ is Willmore if and only if it  satisfies the fourth-order equation $LH+\left(\frac{H^{3}}{2}+2H\right)=0$;}
\item{if $(\Sigma,\g)=(S_{p,\rho},\g_{0})$ is a totally umbilic sphere with the corresponding pullback metric, then the second variation formula is given by 
\be
\label{interm}
\d^{2} I_{0}\left(\Si\right)\left[u_{1},u_{2}\right]=\left(I''_{0}u_{1},u_{2}\right)_{L^{2}(\Si,\g)}, \ \textrm{for} \ \ I''_{0}\left(\Sigma\right)\left[u\right]=\frac{1}{2}\D_{\Si,\g_{0}}\left(\D_{\Si,\g_{0}}+\frac{H^{2}}{2}+2\right)u
\ee
where $H=\frac{\sin(2\rho)}{\sin^{2}(\rho)}$.}
\end{enumerate}

\end{lem}

\begin{rem}
\label{pullback}
For the purpose of the present work, it is  convenient to pull back $u$ to the standard unit sphere $\mathbb{S}^{2}$ (notice that in \eqref{interm} the operator $\D_{\Si,\g}$ depends on the metric $\gamma$ induced on the CMC 2-sphere $\Sigma$, while it would be much more convenient to work with a \textsl{normalized operator}, making the dependence on $\rho$ explicit) and to this aim,  we define the following correspondence: 
\be
\label{corr}
u\in C^{4,\a}\left(S_{p,\rho};\R\right)\ \ \leadsto \ \ w\in C^{4,\a}\left(\mathbb{S}^{2};\R\right) \ : \ w\left(\Theta\right)=u\left(\textrm{exp}_{p}\left(\rho\Theta\right)\right).  
\ee 
For the sake of clarity, let us set $f=\textrm{exp}_{p}\left(\rho\cdot\right)$ so that one simply has $w=u\circ f$.
Then by the scaling properties of the Laplace-Beltrami operator (and the Gauss Lemma) we get
\be\nonumber
\D_{\Si,\g}u\left(q\right)=\frac{1}{\sin^{2}\left(\rho\right)}\D_{\mathbb{S}^{2}}w\left(f^{-1}\left(q\right)\right)
\ee 
and so, as a result, our second derivative operator takes the final form
\be\nonumber
I''_{0}\left[w\right]=\frac{1}{2\sin^{4}\left(\rho\right)}\D_{\mathbb{S}^{2}}^{2}w+\frac{1}{\sin^{2}\left(\rho\right)}\D_{\mathbb{S}^{2}}w+\frac{\sin^{2}\left(2\rho\right)}{4\sin^{6}\left(\rho\right)}\D_{\mathbb{S}^{2}}w;
\ee
by the well-known trigonometric identity $\sin\left(2\rho\right)=2\sin\left(\rho\right)\cos\left(\rho\right)$, we end up getting
\be
\label{secvar}
I''_{0}\left(S_{p,\rho}\right)\left[w\right]=\frac{1}{2\sin^{4}\left(\rho\right)}\D_{\mathbb{S}^{2}}\left(\D_{\mathbb{S}^{2}}+2\right)w.
\ee
Notice that here we are identifying the spaces $C^{4,\a}\left(S_{p,\rho};\R\right)$ and $C^{4,\a}\left(\mathbb{S}^{2};\R\right)$, so that the functional $I_{0}$ is in fact defined on the latter of these (coherently with \cite{Mon1},\cite{Mon2}) and we will always stick to this convention in the sequel.
As a further remark, observe that the operator $I''_{0}\left[w\right]$ has all the scaling and symmetry properties we might expect and, more specifically, it is invariant under the map $\rho\mapsto \pi-\rho$, as it must be.
\end{rem}

\begin{rem}
\label{nucleo}
It is easily checked from \eqref{secvar} that the operator $I''_{0}(S_{p,\rho})$ is  Fredholm of index 0; moreover its kernel $K\subset L^2(\Ss^2,\R)$  is given by  the linear span $\left\langle 1, x_{1}, x_{2}, x_{3}\right\rangle$, where  $x_{1}, x_{2}, x_{3}$ are the restrictions of the coordinate functions of $\R^{3}$ to $\mathbb{S}^{2}\hookrightarrow \R^{3}$. 
\end{rem}

\subsection{Perturbation methods: the Lyapunov-Schmidt reduction}
\label{LS}

The most basic idea behind our approach is to find critical points of the functional $I_{\e}$ by applying a \textsl{finite dimensional reduction}, after which our main theorem will follow by showing that a certain function of four variables defined on $Z$ has an interior maximum point. The tool we need is a sort of \textsl{generalized implicit function theorem}, which is usually referred to as \textsl{Lyapunov-Schmidt reduction}. We recall here its general formulation (see \cite{ab1}-\cite{ab2}, and \cite{am2} for a wider discussion of the method).
\

\
Let $H$ be a Hilbert space and let us consider a suitably smooth functional $J_{\e}: H\to\mathbb{R}$ of the form
\be\nonumber
J_{\e}(u)=J_{0}(u)+\e G(u),
\ee
for some $J_{0}\in C^{2}\left(H; \R\right)$ which plays the role of the \textsl{leading term} (namely the unperturbed functional) and where $G\in C^{2}\left(H;\R\right)$ is an additive perturbation. 
Let us \textsl{assume} that $J_{0}$ has a \textsl{finite dimensional} smooth manifold of critical points:
\be\nonumber
\Xi=\left\{\xi\in H \ | \ J_{0}'\left(\xi\right)=0\right\}.
\ee 
The general idea behind the method is that \textsl{if $J_{0}$ satisfies suitable non-degeneracy conditions, then for $\e$ small enough the functional $J_{\e}$ has a finite-dimensional natural constraint, namely there exist a smooth finite dimensional manifold $\Xi_{\e}$ such that the critical point of $J_{\e}$ constrained to $\Xi_{\e}$ are in fact stationary points for $J_{\e}$}.
Such non-degeneracy conditions are:
\begin{enumerate}
\item[i)]{for all $\xi\in\Xi$ one has that $T_{\xi}\Xi=Ker\left(J''_{0}\right)\left(\xi\right)$;}
\item[ii)]{for all $\xi\in\Xi$ the second derivative operator $J''_{0}\left(\xi\right)$ is a Fredholm operator of index 0}
\end{enumerate}
and the precise statement is the following.

\begin{thm} Suppose that the functional $J_{0}$ has a critical manifold of dimension $d$ and satisfies conditions ${\rm i)}$ and  ${\rm ii)}$ above.
Given a compact subset $\Xi_{c}$ of $\Xi$, there exists $\e_{0}>0$ such that for all $\e\in\left(-\e_{0},\e_{0}\right)$ there is a smooth function $w_{\e}:\Xi_{c}\to H$ satisfying the following three properties:
\begin{enumerate}
\item{for $\e=0$ it results $w_{\e}(\xi)=0,$ for all $\xi\in\Xi_{c}$;}
\item{$w_{\e}(\xi)$ is orthogonal to $T_{\xi}\Xi\,$ for all $\xi\in\Xi_{c}$;}
\item{the manifold $\Xi_{\e}=\left\{\xi+w_{\e}(\xi):\ \xi\in\Xi_{c}\right\}$ is a natural constraint for $J_{\e}$, by which we mean that if $\xi_{\e}$ is a critical point for the function $\Phi_{\e}: \Xi_{c}\to\mathbb{R}$ given by $\Phi_{\e}(\xi)=J_{\e}\left(\xi+w_{\e}(\xi)\right)$, then $u_{\e}=\xi_{\e}+w_{\e}\left(\xi_{\e}\right)$ is a critical point of $J_{\e}$.}
\end{enumerate}
\end{thm}

\begin{rem}
When applying this method, it is often difficult to characterize the map $w$, so that it is in fact necessary to upgrade this scheme showing that under suitable regularity assumptions on $J_{\e}$ the function $w_{\e}(\xi)$ is of order $O(\e)$ uniformly for $\xi\in\Xi_{c}$ (possibly depending on some parameters), so that as a result $\Phi(\xi)=J_{\e}(\xi)+o(\e)$ and so it is sufficient to study the perturbed functional $J_{\e}$ on the critical manifold $\Xi$ (which is typically known).
\end{rem}

\section{Finite-dimensional reduction of the problem}
\label{fdred}
In this section we state and prove the two key Lemmas that allow the finite-dimensional reduction of our problem.

\

The goal of this work is to solve, for suitably small $\e>0$ the Willmore equation $I_{\e}'\left(S_{p,\rho}\left(w\right)\right)=0$ in $w\in C^{4,\a}\left(\mathbb{S}^{2};\R\right).$ Our \textsl{ansatz} is that in fact we can split this into two problems, namely
\be
\begin{cases}
PI_{\e}'\left(S_{p,\rho}\left(w\right)\right)=0\\
\left(I_{C^{0,\a}\left(\mathbb{S}^{2};\R\right)}-P\right)I_{\e}'\left(S_{p,\rho}\left(w\right)\right)=0,
\end{cases}
\ee
where $P:C^{0,\a}\left(\mathbb{S}^{2};\R\right) \to C^{0,\a}\left(\mathbb{S}^{2};\R\right)$ is the projection operator defined in Section 2.1, after equation \eqref{split}. In this section we are concerned with the first of the two, called \textsl{auxiliary equation}; as it will be clear in the sequel, this equation is somehow simpler even if  infinite-dimensional.\newline
In the next lemma we show, using an implicit function type argument, that such equation is solvable at least for suitably small perturbative parameter, i.e. $0<\e<<1$. 

\begin{lem}
\label{LSred}
For each suitably small $\d>0$ there exists $\e_{0}=\e_{0}\left(\d\right)$ and $r_{0}>0$ such that, for $\e\in\left(-\e_{0},\e_{0}\right)$ and $\rho\in\left[\d,\pi-\d\right]$, the auxiliary equation $PI_{\e}'\left(S_{p,\rho}\left(w\right)\right)=0$ has a unique solution $w_{\e}\left(p,\rho\right)\in B\left(0,r_{0}\right)\subseteq C^{4,\a}\left(\mathbb{S}^{2};\R\right)^{\perp}$. Moreover,
\begin{enumerate}
\item{the map $w_{\e}\left(\cdot,\cdot\right):\mathbb{S}^{3}\times \left[\d,\pi-\d\right]\to C^{4,\a}\left(\mathbb{S}^{2};\R\right)^{\perp}$ is $C^{1}$;}
\item{$\left\|w_{\e}\left(p,\rho\right)\right\|_{C^{4,\a}\left(\mathbb{S}^{2};\R\right)}=O\left(\e\right)$ uniformly in $\left(p,\rho\right)\in \mathbb{S}^{3}\times \left[\d,\pi-\d\right].$}
\end{enumerate}
\end{lem}

\begin{pfn}
Given the first variation formula \eqref{first}, which has been derived for a generic isometric immersion of $(\Sigma,\gamma)$ into a Riemannian 3-manifold $(M,g)$, it is immediate to notice that in the special case of $M=\mathbb{S}^{3}$ and $g_{\e}=g_0+\e h$, then
one has 
\be
\label{aux1}
I_{\e}'\left(S_{p,\rho}\left(w\right)\right)=I_{0}'\left(S_{p,\rho}\left(w\right)\right)+\e G\left(\e, S_{p,\rho}\left(w\right)\right)
\ee
where $G\left(\cdot,\cdot\right)$ is a smooth function which is uniformly bounded for $\e$ suitably small. It should be remarked that even though $G\left(\e,\cdot\right)$ is defined on the perturbed sphere $S_{p,\rho}(w)\subseteq \mathbb{S}^{3}$, it is convenient (with slight abuse of notation) to consider it defined on $\mathbb{S}^{2}$ instead: this being said, we observe that $G$ is a smooth function of
 $\e, \Theta, w, Dw, D^{2}w, D^{3}w, D^{4}w$ where $D^{\left(k\right)}w$ synthetically denotes the $k$-th order covariant derivatives of the function $w\in C^{4,\a}\left(\mathbb{S}^{2};\R\right)$.
Let us now define the remainder term $R_{p,\rho}\left(w\right)=I_{0}'\left(S_{p,\rho}\left(w\right)\right)-I_{0}''\left(S_{p,\rho}\right)\left[w\right]$ so that, thanks to \eqref{aux1}, the auxiliary equation we aim at solving takes the expanded form
\be
\label{aux2}
PI_{0}''\left(S_{p,\rho}\right)\left[w\right]+PR_{p,\rho}\left(w\right)+\e PG\left(\e, S_{p,\rho}\left(w\right)\right)=0.
\ee
At this point, in order to turn our problem into a fixed point equation (to be solved by iterative schemes), we recall that the second derivative operator $I_{0}''\left(S_{p,\rho}\right)$ given in  \eqref{secvar}, takes the form
\be\nonumber
I_{0}''\left(S_{p,\rho}\right)\left[w\right]=\frac{1}{2\sin^{4}\left(\rho\right)}\D_{\mathbb{S}^{2}}\left(\D_{\mathbb{S}^{2}}+2\right)w.
\ee
Observe that, if $f \in C^{0,\alpha}(\Ss^2,\R)^\perp$ and  $u \in C^{4,\alpha}(\Ss^2,\R)^\perp$ solves $I_0''(S_{p,\rho})[u]=f$, then (looking at $I_0''$ as a composition of two \textsl{linear} bounded second-order elliptic operators and applying the  Schauder estimates)
\be
\label{schauder}
\left\|u\right\|_{C^{4,\a}\left(\Ss^{2};\R\right)}\leq C \left\|f\right\|_{C^{0,\a}\left(\Ss^{2};\R\right)},
\ee
for some constant $C=C\left(\a,\delta\right)$. Therefore $I_{0}''\left(S_{p,\rho}\right)^{-1}:C^{0,\a}\left(\Ss^{2};\R\right)^{\perp}\to C^{4,\a}\left(\Ss^{2};\R\right)^{\perp}$ is well-defined as a bounded operator and  we can turn \eqref{aux2} into the equivalent \textsl{fixed point problem}
\be\nonumber
w=F_{\e,p,\rho}\left(w\right), \ \textrm{where} \ F_{\e,p,\rho}\left(w\right)=-I_{0}''\left(S_{p,\rho}\right)^{-1}\left[\e PG\left(\e,S_{p,\rho}\left(w\right)\right)+PR_{p,\rho}\left(w\right)\right].
\ee 
Henceforth the only issue is to show that, given $\d$ as in our statement, we can find a threshold $\e_{0}$ so that for any $\left|\e\right|<\e_{0}$, $p\in \mathbb{S}^{3}$, and $\rho\in\left[\d,\pi-\d\right]$, the map $F_{\e,p,\rho}: B\left(0,r\right)\subseteq C^{4,\a}\left(\mathbb{S}^{2};\R\right) \to B(0,r)$ is a contraction (for some $r>0$ small, to be determined).

Now, it is convenient to proceed in two steps:
\begin{enumerate}
\item[\textsl{Step 1:}]{we show that there exists positive $\e$ and $r$ so small that $F_{\e,p,\rho}\left(w\right)$ maps the ball $B\left(0,r\right)\subseteq C^{4,\a}\left(\mathbb{S}^{2};\R\right)$ into itself}
\item[\textsl{Step 2:}]{we show that, possibly by further decreasing $\e$ and $r$ with respect to step 1, we have that in fact that $F_{\e,p,\rho}\left(w\right)$ acts as a contraction on such ball.}
\end{enumerate}

Concerning \textsl{Step 1}, we can make use of the argument given in \cite{Mon1} (pg. 605-606) to show that given $\d$ as above there exists a constant $C(\d)$ such that
\be\nonumber
\left\|F_{\e,p,\rho}\left(w\right)\right\|_{C^{4,\a}\left(\Ss^{2};\R\right)}\leq C\left(\d\right)\left\|PI''_{0}\left(S_{p,\rho}\right)^{-1}\right\|\left(\e+\left\|w\right\|^{2}_{C^{4,\a}\left(\Ss^{2};\R\right)}\right)
\ee
and therefore, given the fact that the right-hand side of the previous inequality is quadratic in $w$, we can certainly pick $r$ and $\e$ so that $F_{\e,p,\rho}$ is a self-mapping of $B\left(0,r\right)$. \newline
Concerning \textsl{Step 2}, we can make use of equation (38) in \cite{Mon1} (its simple proof can be repeated verbatim) to get
\be\nonumber
\left\|F_{\e,p,\rho}\left(w_{2}\right)-F_{\e,p,\rho}\left(w_{1}\right)\right\|_{C^{4,\a}\left(\Ss^{2};\R\right)}\leq 
C\left(\d\right) \left\|PI''_{0}\left(S_{p,\rho}\right)^{-1}\right\|\left(\e+2r\right)\left\|w_{2}-w_{1}\right\|_{C^{4,\a}\left(\Ss^{2};\R\right)}
\ee
which implies the claim.
As a result, we have shown that there exists $\e_{0}$ such that for $\left|\e\right|<\e_{0}$ the auxiliary equation $PI_{\e}'\left(S_{p,\rho}\left(w\right)\right)=0$ is uniquely solvable in $w$. Now, thanks to the well-known version of the Contraction Mapping Theorem in dependence of parameters (specifically: \cite{bressan} (pp. 22-23) and \cite{amn} (pp. 447-449)) one proves that $w_{\e}\left(p,\rho\right)$ is $C^{0}$ in $\left(\e,p,\rho\right)$ and $C^{1}$ in $\left(p,\rho\right)$. Clearly, this implies (just by uniform continuity) that
\be\nonumber
\lim_{\e\to 0}\left\|w_{\e}\left(p,\rho\right)\right\|_{C^{4,\a}\left(\Ss^{2};\R\right)}=0, \ \textrm{uniformly for} \ p\in \mathbb{S}^{3}, \rho\in\left[\d,\pi-\d\right].
\ee
We need to improve this result to show that in fact $\left\|w_{\e}\left(p,\rho\right)\right\|_{C^{4,\a}\left(\Ss^{2};\R\right)}=O(\e)$ uniformly on our domain, which would end the proof. To this aim, let us recall from \eqref{aux2} that the auxiliary equation takes the form
\be\nonumber
PI_{0}''\left(S_{p,\rho}\right)\left[w\right]+\frac{1}{\sin^{4}\left(\rho\right)}Q^{\left(2\right)\left(4\right)}\left(w\right)=-\e PG\left(\e, S_{p,\rho}\left(w\right)\right),
\ee
which is more conveniently analyzed after dividing by $\e$ namely in the form
\be\nonumber
PI_{0}''\left(S_{p,\rho}\right)\left[\frac{w}{\e}\right]+\frac{1}{\e\sin^{4}\left(\rho\right)}Q^{\left(2\right)\left(4\right)}\left(w\right)=- PG\left(\e, S_{p,\rho}\left(w\right)\right).
\ee
Indeed, the right-hand side is uniformly bounded in $C^{0,\a}$ norm for $\left|\e\right|<\e_{0}$ and so the left-hand side has to be as well, but clearly the second summand is of order strictly higher than the first (since obviously $\left\|Q^{\left(2\right)\left(4\right)}\left(w\right)\right\|_{C^{0,\a}\left(\mathbb{S}^{2};\R\right)}\leq C\left\|w\right\|^{2}_{C^{4,\a}\left(\mathbb,{S}^{2};\R\right)}$) and, as a result, $PI_{0}''\left(S_{p,\rho}\right)\left[\frac{w}{\e}\right]$ has to be uniformly bounded in $C^{0,\a}$ for $\e\to 0$. Hence, it follows from Schauder estimates that the function $\frac{w}{\e}$ is bounded in $C^{4,\a}-$norm, which is equivalent to $\left\|w_{\e}\left(p,\rho\right)\right\|_{C^{4,\a}\left(\Ss^{2};\R\right)}=O(\e)$, as we claimed.
\end{pfn}

Following the general Lyapunov-Schmidt reduction as our model (as outlined in Section \ref{LS}), we will now show that Lemma \ref{LSred} determines a \textsl{natural constraint}, in the sense that the problem of proving existence of (conformal) Willmore surfaces is reduced to finding critical points of the $C^{1}$-function $\Phi_{\e}:\mathbb{S}^{3}\times \left[\d,\pi-\d\right]\to\R$ given by
\be
\label{funzrid}
\Phi_{\e}(p,\rho):=I_{\e}\left(S_{p,\rho}\left(w_{\e}\left(p,\rho\right)\right)\right).
\ee

\begin{lem}
\label{natur}
Given $\d>0$ suitably small, let $\e_{0}$ and $r_{0}$ be given as in the statement of Lemma \ref{LSred} and let $\Phi_{\e}$, for $\e\in\left(-\e_{0},\e_{0}\right)$ be the reduced functional defined above by equation \eqref{funzrid}.
Then there exists $\e_{0}'\in\left(0,\e_{0}\right)$ such that if $|\e|<\e_{0}'$ and $\Phi_{\e}$ has a critical point $\left(p_{\e},\rho_{\e}\right)\in \mathbb{S}^{3}\times\left(\d,\pi-\d\right)$, then $S_{p_{\e},\rho_{\e}}\left(w_{\e}\left(p_{\e},\rho_{\e}\right)\right)$ is a critical point for $I_{\e}$.
\end{lem}

\begin{pfn}
By construction, we already know that $S_{p_{\e},\rho_{\e}}\left(w_{\e}\left(p_{\e},\rho_{\e}\right)\right)$ solves the auxiliary equation 
\be\nonumber
PI'_{\e}\left(S_{p_{\e},\rho_{\e}}\left(w_{\e}\left(p_{\e},\rho_{\e}\right)\right)\right)=0
\ee
and so we only need to show that the orthogonal component of $I'_{\e}\left(S_{p_{\e},\rho_{\e}}\left(w_{\e}\left(p_{\e},\rho_{\e}\right)\right)\right)\in L^{2}\left(\mathbb{S}^{2}\right)$ vanishes as well, namely that $\left(I_{L^2\left(\mathbb{S}^{2};\R\right)}-P\right)I_{\e}'\left(S_{p,\rho}\left(w\right)\right)=0$. To this aim, let us recall from remark \eqref{nucleo} (based on the explicit formula \eqref{secvar}) that $K=\textrm{Ker}\left[\D_{\mathbb{S}^{2}}\left(\D_{\mathbb{S}^{2}}+2\right)\right]$ is spanned over $\R$ by (the restriction to $\mathbb{S}^{2}\subset\R^{3}$ of) constant and affine functions; let us consider an orthonormal basis $\{q_{i}^{\e}\}_{i=1,\ldots,4}$ for $K$ with respect to the $L^{2}$-inner product, obtained by normalizing such functions:
\be\nonumber
q_{0}=\frac{1}{\sqrt{4\pi}}, \ q_{i}=\frac{1}{\sqrt{4\pi}}\frac{x_{i}}{|x_{i}|} \ \textrm{for} \ i=1,2,3.
\ee
Let us then decompose $I'_{\e}\left(S_{p_{\e},\rho_{\e}}\left(w_{\e}\left(p_{\e},\rho_{\e}\right)\right)\right)$ with respect to  this basis:
\be\nonumber
I'_{\e}\left(S_{p_{\e},\rho_{\e}}\left(w_{\e}\left(p_{\e},\rho_{\e}\right)\right)\right)=\sum_{i=1}^{4}A_{i,\e}q_{i}, \ A_{i,\e}=\left(I'_{\e}\left(S_{p_{\e},\rho_{\e}}\left(w_{\e}\left(p_{\e},\rho_{\e}\right)\right)\right), q_{i}\right)_{L^{2}\left(\mathbb{S}^{2}\right)}
\ee
so that the assertion we need to prove reduces to showing that $A_{i,\e}=0$ for $i=0,1,2,3.$
To that aim let us make explicit the condition that $\left(p_{\e},\rho_{\e}\right)$ is a stationary point for $\Phi_{\e}$:
\be\nonumber
\left[\frac{\partial \Phi_{\e}}{\partial \rho}\right]\left(p_{\e},\rho_{\e}\right)=0, \ \frac{\partial \Phi_{\e}}{\partial p_{i}}\left(p_{\e},\rho_{\e}\right)=0,
\ee  
where $\frac{\partial}{\partial p_{i}}$ is computed in local coordinates around $p_{\e}$ (and, specifically, by taking geodesic normal coordinates centered at $p_{\e}$ so that $\mathbb{S}^{3}\setminus{\hat{p}_{\e}}$ is identified with $\R^{3}$). For brevity, let us refer from now onwards to the variable $\rho$ as $p_{0}$ so that $\frac{\partial}{\partial p_{0}}$ will stand for $\frac{\partial}{\partial\rho}$ (this allows to use a unified notation). Therefore we have:
\be\nonumber
0=dI_{\e}\left(S_{p_{\e},\rho_{\e}}\left(w_{\e}\left(p_{\e},\rho_{\e}\right)\right)\right)\left[X^{\left(i\right)}\right], \ \textrm{for} \ i=0,1,2,3,
\ee
where $X^{\left(i\right)}$ is the variation vector field (with respect to $p_{i}$) of $\left(S_{p_{\e},\rho_{\e}}\left(w_{\e}\left(p_{\e},\rho_{\e}\right)\right)\right)$. We can decompose each of these vector fields (which are sections, defined over $S_{p_{\e},\rho_{\e}}\left(w_{\e}\left(p_{\e},\rho_{\e}\right)\right)$ of the tangent bundle $T\mathbb{S}^{3}$) into their tangential and normal component: namely, if $\nu$ is an (outward-pointing) co-normal vector field along  $S_{p_{\e},\rho_{\e}}\left(w_{\e}\left(p_{\e},\rho_{\e}\right)\right)$ and $\tau_{1},\tau_{2}$ are a local orthonormal basis to the tangent space of the same sphere, we can write $X^{\left(i\right)}=X^{\left(i\right)}_{1}\tau_{1}+X^{\left(i\right)}_{2}\tau_{2}+X^{i}_{n}\nu$.
Now clearly $dI_{\e}\left(S_{p_{\e},\rho_{\e}}\left(w_{\e}\left(p_{\e},\rho_{\e}\right)\right)\right)\left[ X^{\left(i\right)}_{l}\tau_{l}\right]=0, \ \textrm{for} \ i=0,1,2,3 \ \textrm{and} \ l=1,2$ because the tangential components of $X^{(i)}$ only determine a re-parametrization of the same sphere $S_{p_{\e},\rho_{\e}}\left(w_{\e}\left(p_{\e},\rho_{\e}\right)\right)$. Therefore, we can reduce our analysis to the normal component, for which we have
\be
\label{critic}
\left(I_{\e}'\left(S_{p_{\e},\rho_{\e}}\left(w_{\e}\left(p_{\e},\rho_{\e}\right)\right)\right), X^{\left(i\right)}_{n}\right)_{L^{2}\left(S_{p_{\e},\rho_{\e}}\left(w_{\e}\left(p_{\e},\rho_{\e}\right)\right)\right)}=0, \ \textrm{for} \ i=0,1,2,3.
\ee 
At this stage, using the fact that $\left\|w_{\e}\left(p,\rho\right)\right\|_{C^{4,\a}\left(\mathbb{S}^{2}\right)}=O_{\rho}(\e)$ for $\e\to 0$ (which was proved in Lemma \ref{LSred}) and recalling equation \eqref{graph}, namely $\Psi_{\e,p,\rho,w}\left(\Theta\right)=\textrm{exp}_{p}\left(\left(\rho+w\left(\Theta\right)\Theta\right)\nota \right) \fn$, we get that
\be
\label{norc}
X^{\left(i\right)}_{n}\left(\Psi_{p,\rho,w, \e}\left(\Theta\right)\right)=\sqrt{4\pi}q_{i}+\frac{\partial}{\partial p_{i}}\left(w_{\e}\left(p,\rho\right)\left(\Theta\right)\right)+O_{\rho}\left(\e\right), i=0,1,2,3
\ee
with the remainder term uniformly bounded for $\rho\in\left[\d,\pi-\d\right]$, as it is in our case.
Hence, going back to \eqref{critic} and using \eqref{norc}, it follows that
\be
\label{sist}
0=\sin^{2}\left(\rho_{\e}\right)\left[\sqrt{4\pi}A_{i,\e}+\sum_{j=0}^{3}A_{j,\e}\left(q_{j},\frac{\partial}{\partial p_{i}}\left(w_{\e}\left(p_{\e},\rho_{\e}\right)\right)\right)_{L^{2}\left(\mathbb{S}^{2}\right)}\right]+O_{\rho}\left(\e\right), \ i=0,1,2,3.
\ee
Since the system above for $A_{i,\e}, i=0,1,2,3$ is \textsl{homogeneous}, the claim will follow from showing that the $4\times 4$ matrix given by $\sin^{2}\left(\rho_{\e}\right)\left[\sqrt{4\pi}\d_{ij}+\left(q_{j},\frac{\partial}{\partial p_{i}}\left(w_{\e}\left(p_{\e},\rho_{\e}\right)\right)\right)_{L^{2}\left(\mathbb{S}^{2}\right)}\right]+O_{\rho}\left(\e\right)$ is non-singular. From the conditions $\left(w_{\e}\left(p_{\e},\rho_{\e}\right), q_{i}\right)_{L^{2}\left(\mathbb{S}^{2}\right)}=0$, by differentiating with respect to $\rho$ and $p_{i}$ for $i=1,2,3$, we get that
\be\nonumber
\left(\frac{\partial}{\partial p_{i}}w_{\e}\left(p_{\e},\rho_{\e}\right),q_{j}\right)_{L^{2}\left(\mathbb{S}^{2}\right)}+\left(w_{\e}\left(p_{\e},\rho_{\e}\right), \frac{\partial q_{j}}{\partial p_{i}}\right)_{L^{2}\left(\mathbb{S}^{2}\right)}=0
\ee 
and hence, recalling once again that $\left\|w_{\e}\left(p_{\e},\rho_{\e}\right)\right\|_{C^{4,\a}\left(\mathbb{S}^{2}\right)}=O_{\rho}\left(\e\right)$ we obtain that 
also $\left(\frac{\partial}{\partial p_{i}}w_{\e}\left(p_{\e},\rho_{\e}\right),q_{j}\right)_{L^{2}\left(\mathbb{S}^{2}\right)}=O_{\rho}\left(\e\right)$.
As a result, possibly taking $\e_{0}'$ smaller than $\e_{0}$ if $\e\in\left(-\e_{0}', \e_{0}'\right)$ the determinant of this matrix is not zero and so system \eqref{sist} forces $A_{i,\e}=0$ for $i=0,1,2,3$ which is what we had to prove.
\end{pfn}

\section{Small radii sharp asymptotics}
\label{srasy}

As should be clear from the proof of Lemma \ref{LSred}, the reduction performed in the previous Section does not extend to the closure of the critical manifold $Z=\mathbb{S}^{3}\times\left[0,\pi\right]$ and therefore (as anticipated in the Introduction), we need a different method to extend it suitably till the bases of such \textsl{cylinder}, where the spheres degenerate to points. To this aim, the strategy is essentially to \textsl{fix} the parameter $\e$ (in a suitably small neighborhood of zero) and to use the radius $\rho$ as perturbative parameter. This construction
follows directly from some results proved in \cite{Mon2} about the expansion of the functional $I$ on suitably small perturbed geodesic spheres in a given Riemannian 3-manifold $(M,g)$ (which happens to be, in our case, $\left(\mathbb{S}^{3},g_{\e}=g+\e h\right)$). The real issue is then to show that the reduction corresponding map $w_{\e}\left(p,\rho\right)$ is sufficiently smooth in all its parameters (and, specifically, in $\e$), since this is needed to give an upper bound (involving both $\e$ and $\rho$ together) for $I_{\e}$ on such perturbed spheres. Notice that, due to the \textsl{uniqueness} part in the statement of both Lemma \ref{smallred} and Lemma \ref{LSred}, the two construction can be glued together and give a global reduction map. 
Concerning the functional analytic setting, and specifically concerning the definition of the subspace $C^{4,\a}\left(\Ss^{2};\R\right)^{\perp}$ the reader is referred to Section \ref{nota}.

\begin{lem}
\label{smallred}
There exist $\e_{1}, \rho_{1}>0$, $r_{1}>0$ and a $C^{1}$ map $w_{\e}\left(\cdot,\cdot\right):\left[-\e_{1},\e_{1}\right]\times  \mathbb{S}^{3}\times \left[0,\rho_{1}\right]\to C^{4,\a}\left(\Ss^{2};\R\right)^{\perp},$ $\left(\e,p,\rho,\right)\mapsto w_{\e}\left(p,\rho\right)$ such that if $S_{p,\rho}(w)$ is a critical point of the conformal Willmore functional $I_{\e}$ (for some $\e\in \left(-\e_{1},\e_{1}\right)$) with $\left(p,\rho,w\right)\in \mathbb{S}^{3}\times \left[0,\rho_{1}\right] \times B(0,r_{1})$ then $w=w_{\e}\left(p,\rho\right)$. Moreover, the following properties are satisfied:
\begin{enumerate}
\item{for any $p\in \mathbb{S}^{3}$, the map $\left(\e,\rho\right)\mapsto w_{\e}\left(p,\rho\right)$ is $C^{\infty}$;}
\item{$\left\|w_{\e}\left(p,\rho\right)\right\|_{C^{4,\a}\left(\Ss^{2};\R\right)}=O_{\e}\left(\rho^{3}\right)$ as $\rho\to 0$ uniformly for $p\in \mathbb{S}^{3}$;}
\item{$\left\|\frac{\partial}{\partial\rho}w_{\e}\left(p,\rho\right)\right\|_{L^{2}\left(\Ss^{2};\R\right)}=O_{\e}\left(\rho^{2}\right)$ uniformly for $p\in \mathbb{S}^{3}$;}
\item{one has that
\be\nonumber
\left\|w_{\e}\left(p,\rho\right)+\left(-\frac{1}{12}\rho^{3}Ric_{p}\left(\Theta,\Theta\right)+\frac{1}{36}\rho^{3}R(p)\right)\right\|_{C^{4,\a}\left(\mathbb{S}^{2};R\right)}=O_{\e}\left(\rho^{4}\right), \ \textrm{as} \ \rho\to 0.
\ee
}
\item{one has that $\left\|w_\e \right\|_{C^{4,\a}\left(\mathbb{S}^{2}\right)}=O(\e)$ uniformly for $\left(p,\rho\right)\in S^{3}\times \left[0,\rho_{1}\right]$, namely there exists a constant $C=C\left(g_{0}, \e_{1}, \rho_{1} \right)$ such that $\left\|w_\e \right\|_{C^{4,\a}\left(\mathbb{S}^{2}\right)}\leq C\e$.}
\end{enumerate}
\end{lem}

\begin{pfn}
The construction of this map was performed, for a \textsl{fixed} Riemannian metric $g$ in \cite{Mon2} (Lemma 3.10). Here we need to show that the map is smooth in the couple $(\e,\rho)$ and that estimate 5 holds. Concerning the first assertion, let us recall from Proposition 3.9 in \cite{Mon2} the explicit expansion of the first derivative of the conformal Willmore functional (the same result could be easily deduced from Section \ref{fs}):
\be
\label{prima}
I'_{\e}\left(S_{p,\rho}\left(w\right)\right)=\frac{1}{2\rho^{4}}\D_{\mathbb{S}^{2}}\left(\D_{\mathbb{S}^{2}}+2\right)w+\frac{1}{\rho^{2}}\left[O_{p,\e}\left(\rho^{0}\right)+\frac{1}{\rho}L_{p}^{\left(4\right)}\left(w\right)+\frac{1}{\rho^{3}}Q_{p}^{\left(2\right)\left(4\right)}\left(w\right)\right].
\ee
It is easily checked from the (relatively straightforward) computation leading to \eqref{prima}, that all terms on the right-hand side are smooth (i.e. $C^{\infty}$) in the triple $\left(\e,p,\rho\right)$. Let us then \textsl{fix} a point $p\in \mathbb{S}^{3}$ and consider the map $F_{p}$ given by
\be\nonumber
F_{p}:\left[-\e_{1},\e_{1}\right]\times\left[0,\rho_{1}\right]\times C^{4,\a}\left(\mathbb{S}^{2}\right)^{\perp}\to C^{0,\a}\left(\mathbb{S}^{2}\right)^{\perp}, \  F_{p}\left(\e,\rho, w\right)=\left(\rho^{4}PI'_{\e}\left(S_{p,\rho}\left(w\right)\right)\right)
\ee
so that we want $w=w_{\e}\left(p,\rho\right)$ to be the function implicitly defined by the equation $F_{p}\left(\e,\rho,w\right)=0$.
Indeed, for $\e=\rho=0$ we get that (directly from \eqref{prima}, multiplying by $\rho^{4}$ and projecting through $P$) $F_{p}\left(0,0,0\right)=0$ and moreover if we take the derivative $\frac{\partial F_{p}}{\partial w}\left(0,0,0\right)$ (in the appropriate sense of Banach Calculus) we obtain
\be\nonumber
\frac{\partial F_{p}}{\partial w}\left(0,0,0\right)=\frac{1}{2}\D_{\mathbb{S}^{2}}\left(\D_{\mathbb{S}^{2}}+2\right)
\ee
which is  invertible from $C^{4,\a}\left(\mathbb{S}^{2}\right)^{\perp}$ to $C^{0,\a}\left(\mathbb{S}^{2}\right)^{\perp}$. As a result, we get that (for that fixed point $p\in \mathbb{S}^{3}$) there exist positive constants $\e_{p}, \rho_{p}$ so that $\left(\e,\rho\right)\to w_{\e}\left(p,\rho\right)$ is smooth, for $\e\in\left(-\e_{p},\e_{p}\right)$ and $\rho\in\left(-\rho_{p},\rho_{p}\right)$. Thanks to the Implicit Function Theorem for functional depending on parameters (see for instance \cite{bressan}) and the compactness of $\mathbb{S}^{3}$ we get a global $C^{1}$ map, which coincides,  (thanks to the local uniqueness property in both constructions) with the map defined in Lemma 3.10 of \cite{Mon2}. Assertions 2,3 and 4 follows then directly from Lemma 3.10 in \cite{Mon2} (notice the difference in the notation of $w$: what we call here $w$ was called $-\rho w$ in \cite{Mon2}; this explains the apparently different statements).

Concerning estimate 5, we can argue as follows. Given $p\in \mathbb{S}^{3}$ and $\rho\in\left[0,\rho_{1}\right]$ let us consider the first-order Taylor expansion 
\be\nonumber
w_{\e}\left(p,\rho\right)=w_{0}\left(p,\rho\right)+\e \frac{\partial w_{\e}}{\partial \e}_{\e=\xi}\left(p,\rho\right)
\ee
and hence, again by a local uniqueness argument we must have $w_{0}\left(p,\rho\right)=0$ (since totally umbilic spheres are trivially critical points for the conformal Willmore functional).
It follows that
\be\nonumber
\left\|w_{\e}\left(p,\rho\right)\right\|_{C^{4,\a}\left(\mathbb{S}^{2}\right)}\leq \e \max_{\e\in\left[-\e_{1},\e_{1}\right]}\left\|\frac{\partial w_{\e}\left(p,\rho\right)}{\partial \e}\right\|_{C^{4,\a}\left(\mathbb{S}^{2}\right)},
\ee
but now it is enough to observe that the composite map $\left[-\e_{1},\e_{1}\right]\times\mathbb{S}^{3}\times \left[0,\rho_{1}\right]\to C^{4,\a}\left( \Ss^{2};\R \right)^{\perp}\to\R$ given by $\left\|\frac{\partial w_{\e}\left(p,\rho\right)}{\partial \e}_{\e=0}\right\|_{C^{4,\a}\left(\mathbb{S}^{2}\right)}$ is $C^{0}$ and is defined on a compact space, so that it attains a finite maximum value and this implies the claim.
\end{pfn}

These results being given, we  analyze the asymptotics (in $\rho$, but depending on the parameter $\e$) of the conformal Willmore functional.

\begin{lem}[Proposition 3.11 in \cite{Mon2}, improved]
\label{smallesp} Let $\e_1, \rho_1$ be given by Lemma \ref{smallred} and let  $p\in \mathbb{S}^{3},$ $\rho\in\left[0,\rho_{1}\right]$ and $\e\in [-\e_{0},\e_0]$. For $g_{\e}=g_0+\e h$, the expansion of the conformal Willmore functional on perturbed geodesics spheres $S_{p,\rho}\left(w_{\e}\left(p,\rho\right)\right)$ (determined by the previous Lemma \ref{smallred}) is 
\be\nonumber
I_{\e}\left(S_{p,\rho}\left(w_{\e}\left(p,\rho\right)\right)\right)=\frac{\pi}{5}\left\|\mathring{Ric}_{g_\e}(p)\right\|^{2}\rho^{4}+\Omega\left(\e,\rho\right)
\ee
with
\be\nonumber
\left|\Omega\left(\e,\rho\right)\right|\leq C\e^{2}\rho^{5}, \ \textrm{for} \ \e<\e_{0} \ \textrm{and} \ \rho<\rho_{0}
\ee
for some constant $C\in\R_{>0}$ that can be chosen independently of $\e$ and $\rho$. 
\end{lem}

From subsection \ref{RGP}, we know that  $\left\|\mathring{Ric}_{g_\e}(p)\right\|^{2}=\e^{2}T_{\left(p\right)}^{\left(2\right)}\left(h\right)+o\left(\e^{2}\right)$ (where $T^{\left(2\right)}_{p}\left(h\right)$ is a non-negative quadratic function in the second derivatives of $h$, as we specified above), so we get that
\be
\label{expans}
I_{\e}\left(S_{p,\rho}\left(w_{\e}\left(p,\rho\right)\right)\right)=\frac{\pi}{5}T^{\left(2\right)}_{p}\left(h\right)\e^{2}\rho^{4}+o\left(\e^{2}\right)\rho^{4}+\Omega\left(\e,\rho\right),
\ee
which will be crucial in the sequel of this work.

\

\begin{pfn}
Thanks to  statement 5 of Lemma \ref{smallred}, we obtain  that
\be\label{eq:Ie=e2}
I_{\e}\left(S_{p,\rho}\left(w_{\e}\left(p,\rho\right)\right)\right)=O_{p,\rho}\left(\e^{2}\right) \text{ uniformly for } \rho<\rho_{0}, \ p\in \mathbb{S}^{3}.
\ee
Indeed, let us first Taylor expand in the perturbative parameter $\e$ to get
\be\nonumber
I_{\e}\left(S_{p,\rho}\left(w_{\e}\left(p,\rho\right)\right)\right)=I_{0}\left(S_{p,\rho}\left(w_{\e}\left(p,\rho\right)\right)\right)+\e G_{1}\left(S_{p,\rho}\left(w_{\e}\left(p,\rho\right)\right)\right)+\e^{2}G_{2}\left(S_{p,\rho}\left(w_{\e}\left(p,\rho\right)\right)\right)+o_{p,\rho}\left(\e^{2}\right)
\ee
and then let us expand in $w=w_{\e}\left(p,\rho\right)$ the first two summands above
\begin{eqnarray}
I_{0}\left(S_{p,\rho}\left(w_{\e}\left(p,\rho\right)\right)\right)&=&I_{0}\left(S_{p,\rho}\right)+I_{0}'\left(S_{p,\rho}\right)\left[w_{\e}\left(p,\rho\right)\right]+\frac{1}{2}I_{0}''\left(S_{p,\rho}\right)\left[w_{\e}\left(p,\rho\right),w_{\e}\left(p,\rho\right)\right]+o_{p,\rho}\left(\e^{2}\right) \nonumber \\
G_{1}\left(S_{p,\rho}\left(w_{\e}\left(p,\rho\right)\right)\right)&=&G_{1}\left(S_{p,\rho}\right)+G_{1}'\left(S_{p,\rho}\right)\left[w_{\e}\left(p,\rho\right)\right]+o_{p,\rho}\left(\e\right). \nonumber
\end{eqnarray}
Now, using statement 5 of Lemma \ref{smallred}  repeatedly we see that all terms are (at least) uniformly quadratic in $\e$ apart from $I_{0}\left(S_{p,\rho}\right), \ I_{0}'\left(S_{p,\rho}\right)\left[w_{\e}\left(p,\rho\right)\right]$ and $G_{1}\left(S_{p,\rho}\right)$ which are all exactly zero and so \eqref{eq:Ie=e2} follows ($I_0(S_{p,\rho})$ and $I_0'(S_{p,\rho})$ are clearly null on the totally umbilic spheres,  the computation of $G_1(S_{p,\rho})$ is analogous to the proof of Lemma 4.6 in \cite{Mon2}).
At this point, let us recall from Proposition $3.11$ in  \cite{Mon2} that in fact we already know that
\be
\label{smoothc}
I_{\e}\left(S_{p,\rho}\left(w_{\e}\left(p,\rho\right)\right)\right)=\frac{\pi}{5}T^{\left(2\right)}_{p}\left(h\right)\e^{2}\rho^{4}+o_{p}\left(\e^{2}\right)\rho^{4}+\Omega\left(\e,\rho\right), \ \Omega\left(\e,\rho\right)=O_{p,\e}\left(\rho^{5}\right)
\ee
so that, by comparison with \eqref{eq:Ie=e2} we obtain that
\be\nonumber
\Omega\left(\e,\rho\right)=O_{p,\e}\left(\rho^{5}\right)=O_{p,\rho}(\e^{2})-\frac{\pi}{5}T^{\left(2\right)}_{p}\left(h\right)\e^{2}\rho^{4}-o_{p}\left(\e^{2}\right)\rho^{4}.
\ee
At this point, let us first observe that $\Omega\left(\e,\rho\right)\in C^{\infty}\left(\left(-\e_{1},\e_{1}\right)\times\left(-\rho_{1},\rho_{1}\right)\right)$ because \textsl{all} other terms in \eqref{smoothc} are (the left-hand side as a consequence of the previous Lemma \ref{smallred}, and the term $o_{p}\left(\e\right)$ because it is a remainder term in the expansion of the traceless Ricci tensor at $p$, hence it does not depend on $\rho$, while it depends analytically on $\e$). Therefore, our claim comes from the following elementary Lemma, whose easy proof (based on a Taylor expansion) is omitted.

\begin{lem}
Given $x^{\ast}, y^{\ast}\in\R_{>0}$ let $f\in\left(C^{h+k}\left(\left[-x^{\ast},x^{\ast}\right]\times\left[-y^{\ast},y^{\ast}\right]\right);\R\right)$
such that 
\be\nonumber
f(x,y)=O_{x}\left(\left|y\right|^{k}\right), \ \textrm{for} \ y\to 0
\ee
and
\be\nonumber
f(x,y)=O_{y}\left(\left|x\right|^{h}\right), \ \textrm{for} \ x\to 0.
\ee
Then there exists a constant $M\in\R_{>0}$ such that $\left|f(x,y)\right|\leq M\left|x\right|^{h}\left|y\right|^{k}$ for all $x\in\left[-x^{\ast},x^{\ast}\right]$ and $y\in\left[-y^{\ast},y^{\ast}\right]$.
\end{lem}

\end{pfn}

\section{Proof of Theorem \ref{main}}
In this section, we give a short and direct proof of the main theorem stated in the Introduction which makes use of the various auxiliary tools developed above.

\

\begin{pfn}
Let $\mathring{Ric}_{g_\e}(p)$ be the traceless Ricci tensor of the Riemannian manifold $\left(\mathbb{S}^{3}, g_{\e}=g_0+\e h\right)$ and let us recall from Subsection \ref{RGP} that, due to the analyticity in $\e$ of all curvature tensors, we have a local expansion of the form
\be\nonumber
|\mathring{Ric}_{g_\e}(p)|^{2}=\sum_{k\geq k_{0}}\e^{k}T^{\left(k\right)}_{p}(h)
\ee
for some $k_{0}\geq 2$ depending on $p$. There are three distinct (and disjoint) cases we need to consider: 
\begin{enumerate}
\item[I)]{there exists a point $\overline{p}\in \mathbb{S}^{3}$ such that $k_{0}(\overline{p})=2$, namely $T^{\left(2\right)}_{\overline{p}}(h)\neq 0$;}
\item[II)]{for all points $p\in \mathbb{S}^{3}$ one has $T^{\left(2\right)}_{p}\left(h\right)=0$ but there exists $\overline{p}\in \mathbb{S}^{3}$ such that $T^{\left(k\right)}_{\overline{p}}\left(h\right)\neq 0$, for some $k\geq 3$;}
\item[III)]{the traceless Ricci tensor vanishes identically on $\left(\mathbb{S}^{3}, g_{\e}\right)$.}
\end{enumerate}
We now develop the proof separately for each of these three cases.

\

\underline{\textbf{Case III: the fully degenerate case}}
First of all, notice that for any $h$ as in the statement of the theorem we can find a constant $\e_{III}$ such that for $\e\in\left(-\e_{III},\e_{III}\right)$ the $(0,2)$ tensor $g_{\e}=g_0+\e h$ actually defines a Riemannian metric.
The assumption $\mathring{\textrm{Ric}}_{g_\e} \equiv 0$ can be restated as
\be\nonumber
\textrm{Ric}_{g_\e}=\frac{1}{3}R_{g_\e}g_{\e}
\ee
so the metric $g_\e$ is Einstein and hence, by the Schur Lemma, it has constant scalar curvature. As a result, the Ricci tensor is parallel and $g_{\e}$ has \textsl{constant sectional curvature} (see for instance \cite{Petersen} pp. 38). It follows that $\left(\mathbb{S}^{3}, g_{\e}\right)$ is homothetic (namely: isometric modulo scaling) to $\mathbb{S}^{3}$ and so the result is trivial (in fact, in this case there is a four-dimensional manifold of critical points for the conformal Willmore functional given by the totally umbilic spheres in $\Ss^3$).  
\

\underline{\textbf{Case I: the non-degenerate case}}
\\Let $\e_1,\rho_1$ be given by Lemma \ref{smallred}. Using the assumption $T^{(2)}_{\bar{p}}\neq 0$ together with Lemma \ref{smallesp}, observe that we can choose $\e_{2}\in (0, \e_1]$ and $\rho_2\in (0, \rho_1]$ such that for every $\e \in (-\e_{2},\e_{2}]$ we have 
\be\label{eq:max}
\max_{\mathbb{S}^{3}\times\left(\left[0,\rho_{2}\right]\cup \left[\pi-\rho_{2},\pi\right]\right)}I_{\e}\left(S_{p,\rho}\left(w_{\e}\left(p,\rho\right)\right)\right)<\frac{\max_{\mathbb{S}^{3}\times\left[0,\pi\right]}I_{\e}\left(S_{p,\rho}\left(w_{\e}\left(p,\rho\right)\right)\right)}{2}.
\ee
As a second step, let us set $\e_{I}=\min\left\{\e_{0}',\e_{2}\right\}$ (where $\e_{0}'$ is given by Lemma \ref{natur}, applied with $\d=\rho_{2}$) and $\overline{\rho}=\rho_{2}/2$. Now, let $\Phi_{\e}$ be the corresponding finite-dimensional reduced functional, namely let us set for \textsl{a fixed value} $\e<\overline{\e}$
\be\nonumber
\Phi_\e: \mathbb{S}^{3}\times \left[\overline{\rho}, \pi-\overline{\rho}\right]\to \mathbb{R}, \ \ \Phi_{\e}\left(p,\rho\right)=I_{\e}\left(S_{p,\rho}\left(w_{\e}\left(p,\rho\right)\right)\right).
\ee 
By compactness, $I_{\e}$ has (at least) a maximum point $\left(\overline{p},\overline{r}\right)\in \mathbb{S}^{3}\times \left[\overline{\rho}, \pi-\overline{\rho}\right]$ and, from \eqref{eq:max} and the definition of $\bar{\rho}$ as $\rho_2/2$, this must be an \textsl{interior} maximum point. As a result, by the \textsl{naturality of the constraint} (Lemma \ref{natur}), we conclude that correspondingly $S_{p,\rho}\left(w_{\e}\left(\overline{p},\overline{\rho}\right)\right)$ is a critical point for the functional $I_{\e}$. Namely this graph gives a conformal Willmore surface in the perturbed metric $g_{\e}$.
By construction, such submanifold is a saddle point of $I_{\e}$ of index 4, yet a possibly degenerate one.

\

\underline{\textbf{Case II: the degenerate case}}
\\
This is the most delicate case, and our strategy is  to reduce ourselves to Case III by applying a \textsl{quantitative} version of the (classical) Schur Lemma.
Indeed, let us suppose that the traceless Ricci tensor has no second-order term in the $\e$-expansion at any point and let $k_{0}>2$ be the minimum integer such that $T^{\left(k_{0}\right)}_{p}\left(h\right)\neq 0$ for some $p\in \mathbb{S}^{3}$. Incidentally, observe that $k_0$ must be even and hence, since $\left|\mathring{\textrm{Ric}}_{g_{\e}}\right|^{2}$ is positive definite, at least equal  to $4$. Let us also denote by $\overline{p}$ a point (which will be fixed from now onwards) where $T^{\left(k_{0}\right)}_{\overline{p}}\left(h\right)\neq 0$, so that
\be\label{AssumDeg}
\left\|\mathring{\textrm{Ric}}_{g_\e}\right\|^2=\e^{k_0} \left(T^{(k_0)}(h)+\sum_{j=1}^\infty \e^{j} T^{(k_0+j)}(h) \right) \ \ \textrm{at all points.} 
\ee
 Starting from the identity
\be\nonumber
\textrm{Ric}=\mathring{\textrm{Ric}}+\frac{1}{3}Rg
\ee 
(which we are going to apply for $g=g_{\e}$) and taking the divergence of both left and right-hand side, we get
\be\nonumber
\d(\textrm{Ric})=\d (\mathring{\textrm{Ric}})+\frac{1}{3}dR
\ee
so that, by means of the \textsl{contracted Bianchi identity} $dR=2 \,\d (\textrm{Ric})$, we have
\be
\label{diffscal}
dR=6 \, \d(\mathring{\textrm{Ric}}).
\ee
Now, given any point $q\in \mathbb{S}^{3}$ let us pick a length-minimizing geodesic connecting $\overline{p}$ to $q$ (in the corresponding Riemannian metric $g=g_{\e}$): if we integrate equation \eqref{diffscal} along that path, using our assumption \eqref{AssumDeg}, we obtain
\be
\label{intscal}
R_{\e}(q)=R_{\e}\left(\overline{p}\right)+\e^{k_{0}/2}\Lambda_{\e}(q),
\ee
for some smooth function $\Lambda_{\e}$ on $\mathbb{S}^{3}$. 
At this point, we can exploit \eqref{intscal} to get information on the full curvature tensor of $(\mathbb{S}^{3}, g_{\e})$. Indeed, let us recall the \textsl{Ricci decomposition}
\be\nonumber
\textrm{Riem}=\frac{R}{2n\left(n-1\right)}g\cdot g +\frac{1}{n-2}\left(\textrm{Ric}-\frac{R}{n}g\right)\cdot g +W
\ee
where $\cdot$ stands for the \textsl{Kulkarni-Nomizu product} of two symmetric 2-tensors and $W$ is the Weyl tensor (see for instance \cite{GHL}, pp. 182);
for $n=3$ the Weyl tensor vanishes and therefore the previous reduces to
\be\label{eq:RicciDec}
\textrm{Riem}=\frac{R}{12}g\cdot g +\mathring{\textrm{Ric}}\cdot g.
\ee
Making use of our assumption \eqref{AssumDeg} and its consequence \eqref{intscal} into \eqref{eq:RicciDec}, we get that
\be
\label{resto}
\textrm{Riem}_{g_\e}=\frac{R_{g_\e}\left(\overline{p}\right)}{12}g_{\e}\cdot g_{\e} +\e^{k_{0}/2}\widetilde{\textrm{Riem}_\e}
\ee
for a suitable $(0,4)$ \textsl{curvature-type} tensor $\widetilde{\textrm{Riem}}_\e$. 
 Set $r(\e)=\sqrt{\frac{6}{R_\e(\overline{p})}}$, observe  that $\left(S^{3}, g_{\e}\right)$ is locally isometric (then \textsl{globally isometric} since they are diffeomorphic) to an $\e^{k_{0}/2}$-perturbation of the round sphere of radius $r\left(\e\right)$, namely we can write
\be
\label{claimrel}
g_\e=r^2(\e) g_{0} +\e^{k_{0}/2}\tilde{h}
\ee
for some analytic, symmetric $(0,2)-$tensor $\tilde{h}$. Indeed, it is well-known (see for instance \cite{Will}, pp. 90-92) that given a point $q\in M=\mathbb{S}^{3}$ and denoted by $x^{1},x^{2},x^{3}$ normal coordinates centered at $q$ we can express an analytic metric $g$ as a convergent power series with coefficients only depending on the curvature tensor and its covariant derivatives \textsl{at the point} $q$:
\be
\label{anesp}
g_{rs}(x)=\d_{rs}+\sum_{n=2}^{\infty}\sum_{i_{1},i_{2},\ldots, i_{n}=1}^{3}E^{\left(n\right)}_{r,s,i_{1}i_{2}\ldots i_{n}}\left(\textrm{Riem}\right)x^{i_{1}}x^{i_{2}}\ldots x^{i_{n}}
\ee
where for each $n\geq 2$ the coefficient $E^{(n)}$ consists of a finite number of summands, each one being the evaluation at $q$ of a  term of the form
\be
\nabla^{\a_{i_1}} \textrm{Riem}\ast \nabla^{\a_{i_2}} \textrm{Riem}\ast\ldots \ast \nabla^{\a_{i_{m_n}}} \textrm{Riem}, \ \text{ for some } \a_{i_j}\in\mathbb{N}_{\geq 0} \ \text{ satisfying }  \sum_{j=1}^{m_n} (\a_{i_j}+2) = n. 
\ee
The first terms in such expansion are well-known:
\be\nonumber
g_{rs}=\d_{rs}-\frac{1}{3}\sum_{i,j=1}^{3}R_{irjs}(q)x^{i}x^{j}-\frac{1}{6}\sum_{i,j,k=1}^{3}\nabla_{i}R_{jrks}x^{i}x^{j}x^{k}
\ee
\be\nonumber
+\frac{1}{120}\sum_{i,j,k,l=1}^{3}\left\{-6\nabla_{ij}^{2}R_{krls}+\frac{16}{3}\sum_{f=1}^{3}R_{irjf}R_{kslf}\right\}\left(q\right)x^{i}x^{j}x^{k}x^{l}.
\ee 
Now, if we compute such expansion for a metric homothetic to $g_{0}$ (specifically of the form $r^{2}(\e)g_{0}$) and for our metric $g_{\e}$ (just based on the estimate \eqref{resto}) and compare them, we immediately get that relation  \eqref{claimrel} holds provided we set
\be
\label{defh}
\tilde{h}_{rs}=\sum_{n=2}^{\infty}\sum_{i_{1},i_{2},\ldots, i_{n}=1}^{3}E^{\left(n\right)}_{r,s,i_{1}i_{2}\ldots i_{n}}\left(\widetilde{\textrm{Riem}_{\e}}\right)x^{i_{1}}x^{i_{2}}\ldots x^{i_{n}}
\ee
\be
=-\frac{1}{3}\sum_{i,j=1}^{3}\left(\widetilde{\textrm{Riem}_{\e}}\right)_{irjs}x^{i}x^{j}-\frac{1}{6}\sum_{i,j,k=1}^{3}\nabla^{g_{\e}}_{i}\left(\widetilde{\textrm{Riem}_{\e}}\right)_{jrks}x^{i}x^{j}x^{k}+\ldots.
\ee
By our analyticity assumption this series converges on a ball of suitable radius and hence it is easy to check that it determines a well-defined symmetric $(0,2)$-tensor $\tilde{h}$.  \newline
At this point, let us consider equation \eqref{claimrel}: if we expand in $\e$ both left-hand side ($g_{\e}=g_{0}+\e h$) and right-hand side and compare the two we get that it implies the existence of two real numbers $c, \tilde{c}\in \R$ (independent of $\e$) so that 
\be\nonumber
h=cg_{0}, \ \tilde{h}=\tilde{c}g_{0}.
\ee
In fact $c$ and $-\tilde{c}$ are just the coefficients of order 1 and $k_{0}/2$ respectively in the $\e$-expansion of the function $r^{2}(\e)$.
Clearly, the first of these two relations imply that in fact such $g_{\e}$ should be \textsl{totally degenerate} (in the sense of Case III), which is a contradiction.

\
Therefore the proof of the assertion follows from the arguments we gave for Case I and Case III.
\end{pfn} 

\begin{rem}
Concerning our Remark \ref{expas} at the end of the Introduction, we need to indicate how to modify the argument above in order to treat the more general case when $g_{\e}=g_{0}+h_{\e}$ with $h_{\e}$ analytic in all of its variables. In that case the proof is exactly the same, with the only substantial difference that Case II cannot, in general, be reduced to Case III based on equation \eqref{claimrel}. Instead, after the deduction of \eqref{claimrel} we simply need to observe that the Case II can be reduced to Case I: we can recover all our auxiliary estimates (and, specifically, Lemma \ref{smallesp}) replacing $\e$ by $\e^{k_{0}/2}/r(\e)^{2}$ (recall that $r(\e)\simeq 1$ and $k_0\geq 4$) and we can complete the proof following the very same argument used in the non-degenerate Case I treated above. 
\end{rem}

\appendix
\section{Appendix}
We give here the proof of Lemma \ref{fslem} stated in Section \ref{fs}.

Let us briefly recall the setting: $(\Si,\g)$ is an isometrically immersed surface, and $F:\Si\times(-\s,\s)\to\mathbb{S}^{3}$ is a smooth variation such that $F(\Si,0)=\Si$ and $\frac{\partial F}{\partial s}\left(\Sigma,0\right)=u\nu$, where $\nu$ is the (co-)normal vector field of $\Si$ in $\mathbb{S}^{3}$ and $u\in C^{4,\alpha}\left(\Si \right)$ (concerning the geometric quantities we follow the notations of  \cite{huisk} and \cite{LMS}).

\

\begin{pfn}
Concerning the \textsl{first variation}, we just need to recall the well-known formulas:

\be
\label{fsvar}
\frac{d\mu_{\g}}{\partial s}=uHd\mu_{\g}, \quad \quad \frac{\partial H}{\partial s}=Lu
\ee   
where evaluation at $s=0$ is tacitly assumed and $L$ is the Jacobi operator of $\Si$, namely
\be
Lu=-\D_{\Si,\g}u-\left(Ric\left(\nu,\nu\right)+\left|A\right|^{2}\right)u
\ee
(notice the sign convention, which might be not entirely conventional).
As a result
\be\nonumber
\d I_{0}\left(\Si\right)=\frac{\partial}{\partial s}_{s=0}\int_{\Si_{s}}\left(\frac{H^{2}}{4}+1\right)\,d\mu_{\gamma}=\int_{\Si}\left[\frac{1}{2}H Lu\,d\mu_{\gamma}+\left(\frac{H^{2}}{4}+1\right)uH\right]\,d\mu_{\gamma}
\ee
and hence, integrating by parts
\be\nonumber
\d I_{0}\left(\Si\right)=\int_{\Si}u\left[\frac{1}{2} LH+\left(\frac{H^{3}}{4}+H\right)\right]\,d\mu_{\gamma}.
\ee
Correspondingly, (conformal) Willmore surfaces are defined by (weakly) satisfying the fourth-order equation
\be
\label{first1}
LH+2\left(\frac{H^{3}}{4}+H\right)=0.
\ee 
Concerning the second variation, we have
\be\nonumber
\d^{2}I_{0}\left(\Si\right)=\int_{\Si}\frac{\partial}{\partial s}_{s=0}\left[\frac{1}{2} LH+\left(\frac{H^{3}}{4}+H\right)\right]u\,d\mu_{\gamma}
+\int_{\Si}\left[\frac{1}{2} LH+\left(\frac{H^{3}}{4}+H\right)\right]\left(\frac{\partial u}{\partial s}+H u^{2}\right)_{s=0}\,d\mu_{\gamma}
\ee
which reduces, for any critical point of the conformal Willmore functional $I_{0}$ to
\be\nonumber
\d^{2}I_{0}\left(\Si\right)=\int_{\Si}\frac{\partial}{\partial s}_{s=0}\left[\frac{1}{2} LH+\left(\frac{H^{3}}{4}+H\right)\right]u\,d\mu_{\gamma}.
\ee
In order to proceed further, let us set
\be\nonumber
\mathcal{I}u=\frac{\partial}{\partial s}_{s=0}\left[\frac{1}{2} LH+\left(\frac{H^{3}}{4}+H\right)\right]
\ee
and observe that
\be\label{eq:calI1}
\mathcal{I}(u)=-\frac{1}{2}\left[\frac{\partial}{\partial s}, \Delta_{\Si,\gamma}\right]H +\frac{1}{2}LLu-\frac{1}{2}H\frac{\partial}{\partial s}Ric\left(\nu,\nu\right)-\frac{1}{2}H\frac{\partial \left|A\right|^{2}}{\partial s}+\left(\frac{3}{4}H^{2}+1\right)Lu
\ee
where $\left[T_{1},T_{2}\right]$ denotes the commutator of two (suitably regular) scalar operators.
At this point, we can make use of the computations done in Section 3 of \cite{LMS} namely
\begin{eqnarray}
\frac{\partial}{\partial s}Ric(\nu,\nu)&=&u\nabla_{\nu}Ric(\nu,\nu)-2Ric\left(\nabla u,\nu\right)\label{eq:pf1} \\
\frac{\partial}{\partial s}\left|A\right|^{2}&=&-2u\textrm{tr} A^{3}-2A_{ij}\nabla^{i}\nabla^{j}u-2uA^{ij}T_{ij}\label{eq:pf2}  \\
\left[\frac{\partial}{\partial s},\Delta_{\Si,\g}\right]z&=&Hg\left(\nabla u,\nabla z\right)-ug\left(\nabla z,\nabla H\right)-2A(\nabla u,\nabla z)-2u Ric(\nabla z,\nu)-2u g(A,\nabla^{2}z), \label{eq:pf3}
\end{eqnarray}
where we set $T_{ij}=Ric_{ij}+G(\nu,\nu)\gamma_{ij}$, with $G=Ric-\left(R/2\right) g$ the (ambient) Einstein tensor.

Since we are interested in computing the second variation on totally umbilic spheres $S_{p,\rho}$ in $(\mathbb{S}^{3},g_0)$ we have the following simplifications: $\nabla Ric(\nu,\nu)=0$, $A^\circ_{ij}\equiv 0$, and  the mean curvature is constant (depending on $\rho$). Plugging \eqref{eq:pf1}, \eqref{eq:pf2} and  \eqref{eq:pf3} in \eqref{eq:calI1} we obtain 
\be\nonumber
\mathcal{I}u=\frac{1}{2}LLu+HRic\left(\nabla u,\nu\right)+H\left[u \textrm{tr}A^{3}+g\left(A,\nabla^{2}u\right)+u g\left(A,T\right)\right]+\left(\frac{3}{4}H^{2}+1\right)Lu.
\ee
Writing $A=A^{0}+\frac{1}{2}H\g$ and using the formulas computed at page 14 of \cite{LMS} we finally get
\be\nonumber
\mathcal{I}u=\frac{1}{2}LLu+\left(\frac{H^{2}}{4}+1\right)Lu+H Ric\left(\nu,\nabla u\right).
\ee
Moreover, thanks to the identity $\textrm{div}A^{0}=\frac{1}{2}\nabla H+Ric(\nu,\cdot)^{\sharp}$ it follows that $Ric(\nu,\nabla u)=0$ for any variation $u\in C^{4,\alpha}(S_{p,\rho})$ so that the previous formula simplifies and, as a result, the \textsl{second derivative operator} for $I_{0}$ at $S_{p,\rho}$ is given by \eqref{interm}, namely
\be\label{eq:I0''pf}
I''_{0}\left[u\right]=\frac{1}{2}L^{2}u+\left(\frac{H^{2}}{4}+1\right)Lu.
\ee
Observing that $I''_{0}$ is $L^{2}-$self adjoint we can  write the associated  second variation as $d^{2}I_{0}\left(u_{1},u_{2}\right)=\left(I''_{0}u_{1},u_{2}\right)_{L^{2}(\Si,\g)}$. 

In order to make these formulas totally explicit, we need to compute the terms $Ric\left(\nu,\nu\right)$ and $\left|A\right|^{2}$ for any totally umbilic 2-sphere $S_{p,\rho}$. Clearly, the round metric $g_{0}$ on $\mathbb{S}^{3}$ is Einstein, so that obviously $Ric=\frac{S}{3}g_{0}$ where $S$ is the corresponding scalar curvature, which is exactly equal to 6 for the \textsl{unit sphere}, so that $Ric(\nu,\nu)=2$ at all points of each of the spheres $S_{p,\rho}$. Concerning the other terms, in a principal orthonormal frame at a given point one has
\begin{displaymath}
A =
\left( \begin{array}{ccc}
H/2 & 0 \\
0 & H/2 \\
\end{array} \right)
\end{displaymath}
so that $\left|A\right|^{2}=\frac{H^{2}}{2}$ and the computation reduces to determining the mean curvature $H$ of $S_{p,\rho}$. To that aim, one can then use the first variation formula for the area functional, and easy computations give
\be
\label{cmed}
H=\frac{\sin(2\rho)}{\sin^{2}(\rho)} \quad \text{so that }\left|A\right|^{2}=\frac{\sin^{2}(2\rho)}{2\sin^{4}(\rho)}.
\ee
Notice that for $\rho\to 0^{+}$ one has $H\simeq 2/\rho$ and $\left|A\right|^{2}\simeq 2/\rho^{2}$.

At this point, starting from \eqref{eq:I0''pf} we can perform simple algebraic computations to write the second variation as
\be\nonumber
I''_{0}\left[u\right]=\frac{1}{2}\D_{\Si,\g}^{2}u+\frac{1}{4}H^{2}\D_{\Si,\g}u+\D_{\Si,\g}u=\frac{1}{2}\D_{\Si,\g}\left(\D_{\Si,\g}+\frac{H^{2}}{2}+2\right)u.
\ee
We complete the proof of Lemma \ref{fslem} by replacing  $H$ with its explicit expression for $S_{p,\rho}$ given in \eqref{cmed}.
\end{pfn}

\end{document}